\documentclass[12pt]{article}
\usepackage{amsfonts}
\usepackage{amssymb}
\usepackage{amsthm}
\usepackage{amsmath}
\usepackage{latexsym}
\usepackage{hyperref}
\usepackage[T1]{fontenc}
\usepackage[latin1]{inputenc}
\usepackage{graphicx}

\newtheorem{Theo}{Theorem}[section]
\newtheorem{Lemm}{Lemma}[section]
\newtheorem{Prop}{Proposition}[section]
\newtheorem{Rema}{Remark}[section]

\newtheorem{Coro}{Corollary}[section]
\numberwithin{equation}{section}
\input{tcilatex}
\begin{document}

\title{On the solvability of forward-backward stochastic differential
equations driven by Teugels Martingales}
\date{\ \ }
\author{Dalila Guerdouh, Nabil Khelfallah, Brahim Mezerdi \\
%EndAName
University of Biskra, Laboratory of Applied Mathematics, \\
Po.\ Box 145 Biskra (07000), Algeria\thanks{
E-mail addresses: dalilajijel18@gmail.com, nabilkhelfallah@yahoo.fr,
bmezerdi@yahoo.fr.}}
\date{}
\maketitle

\begin{abstract}
We deal with a class of fully coupled forward-backward stochastic
differential equations (FBSDE for short), driven by Teugels martingales
associated with some Lévy process. Under some assumptions on the derivatives
of the coefficients, we prove the existence and uniqueness of a global
solution on an arbitrarily large time interval. Moreover, we establish \
stability and comparison theorems for the solutions of such equations. Note
that the present work extends known results by Jianfeng Zhang (Discrete
Contin. Dyn. Syst. Ser. B 6 (2006), no. 4, 927--940), proved for FBSDEs
driven by a Brownian motion, to FBSDEs driven by general Lévy processes.

\end{abstract}

\section{Introduction}

Let $\left( L_{t}\right) _{0\leq t\leq T}$ be a $\mathbb{R}-$valued Lévy
process defined on a complete filtered probability space $\left( \Omega ,%
\mathcal{F},\left( \mathcal{F}_{t}\right) _{t\geq 0},P\right) $ satisfying
the usual conditions. Assume that the Lévy measure $\nu \left( dz\right) $
corresponding to the Lévy process $L_{t}$ satisfies:

$\left( i\right) $ $\int\limits_{\mathbb{R}}\left( 1\wedge z^{2}\right) \nu
\left( dz\right) <\infty $,

$\left( ii\right) $ there exist $\alpha >0$ such that for every $\varepsilon
>0,$%
\begin{equation*}
\int\limits_{\left] -\varepsilon ,\varepsilon \right[ ^{c}}e^{\alpha
\left\vert z\right\vert }\nu \left( dz\right) <\infty .
\end{equation*}%
Assumptions $\left( i\right) $ and $\left( ii\right) $ imply in particular
that the random variable $L\left( t\right) $ has moments of all orders. We
also assume that $\mathcal{F}_{t}=\mathcal{F}_{0}\vee \sigma \left(
L_{s},s\leq t\right) \vee \mathcal{N},$ where $\mathcal{G}_{1}\vee \mathcal{G%
}_{2}$ denotes the $\sigma -$field generated by $\mathcal{G}_{1}\cup 
\mathcal{G}_{2}$ and $\mathcal{N}$ is the totality of the $P$-negligible
sets.

The aim of this work is to prove existence and uniqueness of solutions of
the following coupled forward-backward stochastic differential equation
(FBSDE for short)%
\begin{equation}
\left\{ 
\begin{array}{ll}
X_{t}= & X_{0}+\int\limits_{0}^{t}f\left( s,w,X_{s},Y_{s},Z_{s}\right)
ds+\sum\limits_{i=1}^{\infty }\int\limits_{0}^{t}\sigma ^{i}\left(
s,w,X_{s-},Y_{s-}\right) dH_{s}^{i}, \\ 
Y_{t}= & \varphi \left( X_{T}\right) +\int\limits_{t}^{T}g\left(
s,w,X_{s},Y_{s},Z_{s}\right) ds-\sum\limits_{i=1}^{\infty
}\int\limits_{t}^{T}Z_{s}^{i}dH_{s}^{i},%
\end{array}%
\right.  \label{EQ1}
\end{equation}%
where $t\in \left[ 0,T\right] ,$ $H_{t}=\left( H_{t}^{i}\right)
_{i=1}^{\infty }$ are pairwise strongly orthonormal Teugels martingales
associated with the Lévy process $L_{t}$. For any $%
%TCIMACRO{\U{211d} }%
%BeginExpansion
\mathbb{R}
%EndExpansion
$-valued and $\mathcal{F}_{0}$-measurable random vector $X_{0}$, satisfying $%
\mathbb{E}\left\vert X_{0}\right\vert ^{2}$ $<\infty $, we are looking for
an $%
%TCIMACRO{\U{211d} }%
%BeginExpansion
\mathbb{R}
%EndExpansion
\times 
%TCIMACRO{\U{211d} }%
%BeginExpansion
\mathbb{R}
%EndExpansion
\times l\left( 
%TCIMACRO{\U{211d} }%
%BeginExpansion
\mathbb{R}
%EndExpansion
\right) $-valued solution $\left( X_{t},Y_{t},Z_{t}\right) $ on an
arbitrarily fixed large time duration, which is square-integrable and
adapted with respect to the filtration $\mathcal{F}_{t}$ generated by $L_{t}$
and $\mathcal{F}_{0}$ satisfying%
\begin{equation*}
\mathbb{E}\int\limits_{0}^{t}\left( \left\vert X_{t}\right\vert
^{2}+\left\vert Y_{t}\right\vert ^{2}+\left\vert Z_{t}\right\vert
^{2}\right) dt<\infty .
\end{equation*}

The existence and uniqueness of solutions of FBSDEs without the Teugels part
have been widely studied by many authors (see, e.g. \cite{A}, \cite{D}, \cite%
{HP}, \cite{MP}, \cite{PT}, \cite{PW}, and \cite{Y}). The first study of
FBSDEs has been performed by Antonelli \cite{A} in the early 1990s. The
author has used the contraction mapping technique to obtain a local
existence and uniqueness result in a small time interval. Hu and Peng \cite%
{HP} have used a probabilistic method to establish an existence and
uniqueness result, under certain monotonicity conditions, in the case where
the forward and backward components have the same dimension. Then Hamadène 
\cite{H} improved their result by proving it under weaker monotonicity
assumptions. Peng and Wu provided in \cite{PW} more general results by
extending the two above results, without the restriction on the dimensions
of the forward and backward parts.

In spite of\ the large literature devoted to the Brownian case as we have
mentioned above, there are relatively a few results on FBSDEs driven by
Teugels Martingales. To the best of our knowledge, the first paper dealing
with this kind of equations driven by Lévy processes is \cite{PS}, where the
authors have proved the existence and uniqueness via the solution of its
associated partial integro-differential equation (PIDE for short). Then
Baghery et al. \cite{B} proved under some monotonicity assumptions, the
existence and uniqueness of solutions on an arbitrarily fixed large time
duration.

Motivated by the above results and by imposing an assumption on the
derivatives of the coefficients, introduced by Zhang \cite{Z}, we establish
two main results. We shall first prove the existence and uniqueness of the
solution of the FBSDE \ref{EQ1}, without any restriction on the time
duration. The main idea of the proof is to construct the solution on small
intervals, and then extend it piece by piece to the whole interval. In a
second step, we prove stability and comparison theorems for the solutions.
Let us point out that our work extends the results of Jianfeng Zhang
(Discrete Contin. Dyn. Syst. Ser. B 6 (2006), no. 4, 927--940), to FBSDEs
driven by general Lévy processes. We note that much of the technical
difficulties coming from the Teugels martingales are due to the fact that
the quadratic variation $\left[ H^{i},H^{j}\right] $ is not absolutely
continuous, with respect to the Lebesgue measure. To overcome these
difficulties, we use the fact that the predictable quadratic variation
process$\ \left\langle H^{i},H^{j}\right\rangle _{t}$ is equal to $\delta
_{ij}t$ and that $\left[ H^{i},H^{j}\right] _{t}-\left\langle
H^{i},H^{j}\right\rangle _{t}$ is a martingale.

This paper is organized as follows. In Section $2$, we give some
preliminaries and notations about Teugels martingales. In Section $3$, we
give some assumptions and provide our main results. The proofs are provided
in the last section.

\section{Notations and assumptions}

Let us recall briefly the $L^{2}$ theory of Lévy processes as it is
investigated in Nualart-Schoutens \cite{NS}. A convenient basis for
martingale representation is provided by the so-called Teugels martingales.
This means that this family has the predictable representation property.

Denote by $\Delta L_{t}=L_{t}-L_{t_{-}}$ where%
\begin{equation*}
L_{t_{-}}=\lim_{s\rightarrow t,s<t}L_{s},\text{ \ \ }t>0,
\end{equation*}%
and define the power jump processes by%
\begin{equation*}
L_{t}^{\left( i\right) }=\left\{ 
\begin{array}{c}
L_{t}\text{ \ \ \ \ \ \ \ \ \ \ \ \ \ \ \ \ \ if }i=1; \\ 
\sum\limits_{0<s\leq t}\left( \bigtriangleup L_{s}\right) ^{i}\text{ \ \ \ \
\ if }i\geq 2.%
\end{array}%
\right.
\end{equation*}

If we denote 
\begin{equation*}
Y_{t}^{\left( i\right) }=L_{t}^{\left( i\right) }-\mathbb{E}\left[
L_{t}^{\left( i\right) }\right] ,i\geq 1,
\end{equation*}%
with%
\begin{equation*}
\mathbb{E}\left[ L_{t}^{\left( 1\right) }\right] =\mathbb{E}\left[ L_{t}%
\right] =t\mathbb{E}\left[ L_{1}\right] =tm_{1},
\end{equation*}%
and, for $i\geq 2$%
\begin{equation*}
\mathbb{E}\left[ L_{t}^{\left( i\right) }\right] =\mathbb{E}\left[
\sum\limits_{0<s\leq t}\left( \bigtriangleup L\left( s\right) \right) ^{i}%
\right] =t\int_{-\infty }^{\infty }z^{i}\nu \left( dz\right) =tm_{i}.
\end{equation*}

Then the family of Teugels martingales $\left( H_{t}^{i}\right)
_{i=1}^{\infty },$\ is defined by 
\begin{equation*}
H_{t}^{i}=\sum\limits_{j=1}^{j=i}a_{ij}Y_{t}^{\left( j\right) }.
\end{equation*}

The coefficients $a_{ij}$ correspond to the orthonormalization of the
polynomials $1,$ $x,$ $x^{2},$ $...$ with respect to the measure $\mu \left(
dx\right) =x^{2}\nu \left( dx\right) +\delta _{0}\left( dx\right) $. Then $%
\left( H_{t}^{i}\right) _{i=1}^{\infty }$ is a family of strongly orthogonal
martingales such that $\left\langle H^{i},H^{j}\right\rangle _{t}=\delta
_{ij}.t$ and $\left[ H^{i},H^{j}\right] -\left\langle
H^{i},H^{j}\right\rangle _{t}$ is a martingale, see \cite{NS, P}.

The following lemma which gives some useful properties of the Teugels
martingale will be needed in the sequel.

\begin{Lemm}
\label{lem} $i)$ The process $H_{t}^{i}$ can be represented as follows:%
\begin{equation*}
H_{t}^{i}=q_{i-1}\left( 0\right) B_{t}+\int_{\mathbb{R}}p_{i}\left( x\right) 
\tilde{N}\left( t,dx\right)
\end{equation*}%
where $B_{t}$ be a Brownian motion, and $\tilde{N}\left( t,dx\right) $ is
the compensated Poisson random measure that corresponds to the pure jump
part of $L_{t}$ and the polynomials $q_{i-1}\left( 0\right) $ and $%
p_{i}\left( x\right) $ associated to $L_{t}.$

$ii)$ The polynomials $p_{i}$ and $q_{j}$ are linked by the relation:%
\begin{equation*}
\int_{\mathbb{R}}p_{i}\left( x\right) p_{j}\left( x\right) v\left( dx\right)
=\delta _{ij}-q_{i-1}\left( 0\right) q_{j-1}\left( 0\right) .
\end{equation*}
\end{Lemm}

\textbf{Proof. }See \cite{PS}.\hfill $\square $\bigskip

In the rest of this section, we list all the notations that will be
frequently used throughout this work.

$l^{2}:$ the Hilbert space of real-valued sequences $x=\left( x_{n}\right)
_{n\geq 0}$ with norm 
\begin{equation*}
\left\Vert x\right\Vert =\left( \sum\limits_{i=1}^{\infty }x_{i}\right) ^{%
\frac{1}{2}}<\infty .
\end{equation*}

Let us define

$l^{2}\left( \mathbb{R}\right) :$ the space of $\mathbb{R}$-valued process $%
\left\{ f^{i}\right\} _{i\geq 0}$ such that%
\begin{equation*}
\left( \sum\limits_{i=1}^{\infty }\left\Vert f^{i}\right\Vert _{\mathbb{R}%
}^{2}\right) ^{\frac{1}{2}}<\infty .
\end{equation*}

$l_{\mathcal{F}}^{2}\left( 0,T,\mathbb{R}\right) :$ the Banach space of $%
l^{2}\left( \mathbb{R}\right) -$valued $\mathcal{F}_{t}-$predictable
processes such that%
\begin{equation*}
\left( \mathbb{E}\int_{0}^{T}\sum\limits_{i=1}^{\infty }\left\Vert
f^{i}\left( t\right) \right\Vert _{\mathbb{R}}^{2}\right) ^{\frac{1}{2}%
}<\infty .
\end{equation*}

$\mathcal{S}_{\mathcal{F}}^{2}\left( 0,T,\mathbb{R}\right) :$ the Banach
space of $\mathbb{R}-$valued $\mathcal{F}_{t}-$adapted and càdlàg\ processes
such that%
\begin{equation*}
\left( \mathbb{E}\sup_{0\leq t\leq T}\left\vert f\left( t\right) \right\vert
^{2}\right) ^{\frac{1}{2}}<\infty .
\end{equation*}

$L^{2}\left( \Omega ,\mathcal{F},P,\mathbb{R}\right) :$ the Banach space of $%
\mathbb{R}-$valued, square integrable random variables on $\left( \Omega ,%
\mathcal{F},P\right) .$ Here and in what follows, for notational simplicity,
we shall denote%
\begin{equation*}
\begin{array}{ccc}
\int_{0}^{t}\sigma \left( s,w,X_{s-},Y_{s-}\right) dH_{s} & \text{and} & 
\int_{t}^{T}Z_{s}dH_{s}%
\end{array}%
\end{equation*}

instead of%
\begin{equation*}
\begin{array}{ccc}
\sum\limits_{i=1}^{\infty }\int_{0}^{t}\sigma ^{i}\left(
s,w,X_{s-},Y_{s-}\right) dH_{s}^{i} & \text{and} & \sum\limits_{i=1}^{\infty
}\int_{t}^{T}Z_{s}^{i}dH_{s}^{i}%
\end{array}%
\end{equation*}%
respectively, where $Z_{s}=\left\{ Z_{s}^{i}\right\} _{i=1}^{\infty },$ $%
\sigma _{s}=\left\{ \sigma _{s}^{i}\right\} _{i=1}^{\infty },$ $\sigma ^{i}:%
\left[ 0,T\right] \times \Omega \times \mathbb{R}\times \mathbb{R}%
\rightarrow l^{2}\left( \mathbb{R}\right) .$ Further, for the notational
simplicity, we have suppressed $w$ and we will do so below. We also use the
following notation%
\begin{equation*}
M^{2}\left( 0,T\right) =\mathcal{S}_{\mathcal{F}}^{2}\left( 0,T,\mathbb{R}%
\right) \times \mathcal{S}_{\mathcal{F}}^{2}\left( 0,T,\mathbb{R}\right)
\times l_{\mathcal{F}}^{2}\left( 0,T,\mathbb{R}\right) .
\end{equation*}

The following assumptions will be considered in this paper.

\noindent We suppose that the coefficients 
\begin{equation*}
\begin{array}{l}
f:\left[ 0,T\right] \times \Omega \times \mathbb{R}\times \mathbb{R}\times
l^{2}\left( \mathbb{R}\right) \rightarrow \mathbb{R}, \\ 
\sigma :\left[ 0,T\right] \times \Omega \times \mathbb{R}\times \mathbb{R}%
\rightarrow l^{2}\left( \mathbb{R}\right) , \\ 
g:\left[ 0,T\right] \times \Omega \times \mathbb{R}\times \mathbb{R}\times
l^{2}\left( \mathbb{R}\right) \rightarrow \mathbb{R}, \\ 
\varphi :\Omega \times \mathbb{R}\rightarrow \mathbb{R},%
\end{array}%
\end{equation*}%
are progressively measurable, such that:

\begin{enumerate}
\item[$\left( \mathbf{H}_{1}\right) $] There exist $\lambda ,$ $\lambda
_{0}>0,$ such that $\forall t\in \left[ 0,T\right] ,\forall \left(
x,y,z\right) $ and $\left( x^{\prime },y^{\prime },z^{\prime }\right) $%
\textbf{\ }in\textbf{\ }$\mathbf{%
%TCIMACRO{\U{211d} }%
%BeginExpansion
\mathbb{R}
%EndExpansion
\times 
%TCIMACRO{\U{211d} }%
%BeginExpansion
\mathbb{R}
%EndExpansion
\times }$ $l\left( 
%TCIMACRO{\U{211d} }%
%BeginExpansion
\mathbb{R}
%EndExpansion
\right) $ 
\begin{eqnarray*}
\left\vert f\left( t,x,y,z\right) -f\left( t,x^{\prime },y^{\prime
},z^{\prime }\right) \right\vert &\leq &\lambda \left( \left\vert
x-x^{\prime }\right\vert +\left\vert y-y^{\prime }\right\vert +\left\Vert
z-z^{\prime }\right\Vert _{l^{2}\left( \mathbb{R}\right) }\right) , \\
\left\vert \sigma \left( t,x,y\right) -\sigma \left( t,x^{\prime },y^{\prime
}\right) \right\vert ^{2} &\leq &\lambda ^{2}\left( \left\vert x-x^{\prime
}\right\vert ^{2}+\left\vert y-y^{\prime }\right\vert ^{2}\right) , \\
\left\vert g\left( t,x,y,z\right) -g\left( t,x^{\prime },y^{\prime
},z^{\prime }\right) \right\vert &\leq &\lambda \left( \left\vert
x-x^{\prime }\right\vert +\left\vert y-y^{\prime }\right\vert +\left\Vert
z-z^{\prime }\right\Vert _{l^{2}\left( \mathbb{R}\right) }\right) , \\
\left\vert \varphi \left( x\right) -\varphi \left( x^{\prime }\right)
\right\vert &\leq &\lambda _{0}\left( \left\vert x-x^{\prime }\right\vert
\right) .
\end{eqnarray*}

\item[$(\mathbf{H}_{\mathbf{2}})$] The functions $f,g,\sigma ,\varphi $ are
differentiable with respect to $x,$ $y,$ $z$ with uniformly bounded
derivatives such that%
\begin{equation}
\sigma _{y}f_{z}=0\text{ and \ }f_{y}+\sigma _{x}f_{z}+\sigma _{y}g_{z}=0.
\end{equation}
\end{enumerate}

Let us mention that assumption (\textbf{H}$_{2})$ has been introduced bfor
the first time by Zhang \cite{Z} in the case of FBSDEs without jumps.

\section{ The main results}

\subsection{Existence and uniqueness}

The following theorem gives the existence of a solution in a small time
duration.

\begin{Theo}
\label{STD} Suppose that $\left( \mathbf{H}_{\mathbf{1}}\right) $ is
satisfied. Assume further that%
\begin{equation*}
V_{0}^{2}\overset{\bigtriangleup }{=}\mathbb{E}\left\{ \left\vert
X_{0}\right\vert ^{2}+\left\vert \varphi \left( 0\right) \right\vert
^{2}+\int_{0}^{T}\left[ \left\vert f\left( t,0,0,0\right) \right\vert
^{2}+\left\Vert \sigma \left( t,0,0\right) \right\Vert _{l^{2}\left( \mathbb{%
R}\right) }^{2}+\left\vert g\left( t,0,0,0\right) \right\vert ^{2}\right]
dt\right\} <\infty .
\end{equation*}%
Then, for every\ $\mathcal{F}_{0}$-measurable random vector $X_{0},$\ there
exists a constant $\delta $ depending only on $\lambda $ and $\lambda _{0},$
such that for $T\leq \delta ,$ equation $(\ref{EQ1})$ has a unique solution
which belongs to $M^{2}\left( 0,T\right) $.
\end{Theo}

The following proposition gives a priori estimates, which shows in
particular the continuous dependence of the solution upon the data.

\begin{Prop}
\label{Prop1} Under the same assumptions of the Theorem $\ref{STD}$, there
exist $\delta $ and $C_{0}$ depending on $\lambda $ and $\lambda _{0},$ such
that for $T\leq \delta ,$ the following estimates hold true:

\begin{itemize}
\item[i)] 
\begin{equation*}
\left\Vert \Pi \right\Vert =\mathbb{E}\left( \sup\limits_{0\leq t\leq T}%
\left[ \left\vert X_{t}\right\vert ^{2}+\left\vert Y_{t}\right\vert ^{2}%
\right] +\int_{0}^{T}\left\Vert Z_{t}\right\Vert _{l^{2}\left( \mathbb{R}%
\right) }^{2}dt\right) ^{\frac{1}{2}}\leq C_{0}V_{0}.
\end{equation*}

\item[ii)] 
\begin{equation*}
\mathbb{E}\left\{ \sup\limits_{0\leq t\leq T}\left[ \left\vert
X_{t}\right\vert ^{2p}+\left\vert Y_{t}\right\vert ^{2p}\right] +\left(
\int_{0}^{T}\left\Vert Z_{t}\right\Vert _{l^{2}\left( \mathbb{R}\right)
}^{2}dt\right) ^{p}\right\} <\infty
\end{equation*}
\end{itemize}
\end{Prop}

The next Theorem extends the result in Theorem $\ref{STD}$ to arbitrary
large time duration.

\begin{Theo}
\label{LTD} Assume $\left( \mathbf{H}_{\mathbf{1}}\right) $, $\left( \mathbf{%
H}_{\mathbf{2}}\right) $ and $V_{0}^{2}<\infty .$ Then:

$i)$ Equation $(\ref{EQ1})$ has a unique solution $\Pi \in M^{2}\left(
0,T\right) .$

$ii)$ The following estimate holds%
\begin{equation*}
\left\Vert \Pi \right\Vert ^{2}\leq CV_{0}^{2}.
\end{equation*}
\end{Theo}

\subsection{Stability theorem}

The following results state the stability of the solution of FBSDE $(\ref%
{EQ1})$ with respect to the initial condition and the data. This means that
the solution of equation $(\ref{EQ1})$ does not change too much under small
perturbations of the data. In other words, the trajectories which are close
to each other at specific instant should therefore remain close to each
other at all subsequent instants. To state the next theorem and its
corollary, let us consider $\Pi ^{i},i=0,1$ the solutions of $(\ref{EQ1})$
corresponding to $\left( f^{i},\sigma ^{i},g^{i},\varphi ^{i}\right) $. We
shall consider the following notations, $\Delta \Pi \overset{\Delta }{=}\Pi
^{1}-\Pi ^{0}$ and for any function $h\overset{\Delta }{=}f,\sigma
,g,\varphi $, we set $\Delta h\overset{\Delta }{=}h^{1}-h^{0}.$

\begin{Theo}
\label{stability of solution}Assume that $\left( f^{i},\sigma
^{i},g^{i},\varphi ^{i},X_{0}^{i}\right) ,i=0,1,$ satisfy the same
conditions of Theorem \ref{LTD}. Then%
\begin{equation*}
\left\Vert \Delta \Pi \right\Vert ^{2}\leq C\mathbb{E}\left\{ \left\vert
\Delta X_{0}\right\vert ^{2}+\left\vert \Delta \varphi \left(
X_{T}^{1}\right) \right\vert ^{2}+\int_{0}^{T}\left[ \left\vert \Delta
f\right\vert ^{2}+\left\Vert \Delta \sigma \right\Vert _{l^{2}\left( \mathbb{%
R}\right) }^{2}+\left\vert \Delta g\right\vert ^{2}\right] \left( t,\Pi
_{t}^{1}\right) dt\right\} .
\end{equation*}
\end{Theo}

\begin{Coro}
\label{corol2} Suppose that $\left( f^{n},\sigma ^{n},\varphi
^{n},g^{n},X_{0}^{n}\right) ,$ for $n=0,1...$ satisfy the same conditions of
Theorem $\ref{LTD}$. Moreover assume that:

i) $X_{0}^{n}\rightarrow X_{0}^{0}$ in $L^{2}.$

ii) for $h\overset{\Delta }{=}f,\sigma ,\varphi ,g$ $,$ $h^{n}\left( t,\Pi
\right) \rightarrow h^{0}\left( t,\Pi \right) $ as $n\rightarrow \infty .$

$iii)$ $\mathbb{E}\left\{ \left\vert X_{0}^{n}-X_{0}^{0}\right\vert
^{2}+\left\vert \varphi ^{n}-\varphi ^{0}\right\vert ^{2}\left( 0\right)
+\int_{0}^{T}\left[ \left\vert f^{n}-f^{0}\right\vert ^{2}+\left\Vert \sigma
^{n}-\sigma ^{0}\right\Vert _{l^{2}\left( \mathbb{R}\right) }^{2}+\left\vert
g^{n}-g^{0}\right\vert ^{2}\right] \left( t,0,0,0\right) dt\right\}
\rightarrow 0$

Then if $\Pi ^{n\text{ }}(resp.\Pi )$ denotes the solution of $(\ref{EQ1})$
corresponding to $\left( f^{n},\sigma ^{n},\varphi
^{n},g^{n},X_{0}^{n}\right) $ (resp. $\left( f,\sigma ,\varphi
,g,X_{0}^{0}\right) $, we obtain%
\begin{equation*}
\left\Vert \Pi ^{n}-\Pi ^{0}\right\Vert \rightarrow 0\text{ as }%
n\longrightarrow +\infty .
\end{equation*}
\end{Coro}

\subsection{Comparison theorem}

In what follows we provide, under the same assumptions as for the existence
and uniqueness results, another important result, which is the comparison
theorem. Let $\left( X,Y,Z\right) $ be the solution to the following LFBSDE:%
\begin{equation}
\left\{ 
\begin{array}{ll}
X_{t}= & \int_{0}^{t}\left(
a_{s}^{1}X_{s}+b_{s}^{1}Y_{s}+c_{s}^{1}Z_{s}\right) ds+\int_{0}^{t}\left(
a_{s}^{2}X_{s}+b_{s}^{2}Y_{s}\right) dH_{s}, \\ 
Y_{t}= & PX_{T}+\alpha +\int_{t}^{T}\left(
a_{s}^{3}X_{s}+b_{s}^{3}Y_{s}+c_{s}^{3}Z_{s}+\beta _{s}\right)
ds-\int_{t}^{T}Z_{s}dH_{s}.%
\end{array}%
\right.  \label{EQ7}
\end{equation}%
Then we have the following proposition, which is the linear version of the
next theorem.

\begin{Prop}
\label{Prop3} Assume $\left\vert a_{t}^{i}\right\vert ,\left\vert
b_{t}^{i}\right\vert ,\left\vert c_{t}^{i}\right\vert \leq \lambda ,$ $%
\left\vert P\right\vert \leq \lambda _{0}$ and $(\mathbf{H}_{2})$\ holds
true. Assume further that $\alpha \geq 0$ and $\beta _{s}\geq 0.$ Then%
\begin{equation*}
Y_{0}\geq 0.
\end{equation*}
\end{Prop}

Further we have the following general result. Let $\Pi ^{i},i=0,1,$ be the
solution of the following FBSDE:%
\begin{equation}
\left\{ 
\begin{array}{ll}
X_{t}^{i}= & X_{0}+\int_{0}^{t}f\left( s,\Pi _{s}^{i}\right)
ds+\int_{0}^{t}\sigma \left( s,X_{s-}^{i},Y_{s-}^{i}\right) dH_{s}, \\ 
Y_{t}^{i}= & \varphi ^{i}\left( X_{T}^{i}\right) +\int_{t}^{T}g^{i}\left(
s,\Pi _{s}^{i}\right) ds-\int_{t}^{T}Z_{s}^{i}dH_{s},i=0,1%
\end{array}%
\right.  \label{EQ71}
\end{equation}

\begin{Theo}
\label{comparison of solutions} Let $\Pi ^{i},i=0,1,$ be the solutions of
the FBSDEs (\ref{EQ1}). If

$i)$ $\left( f,\sigma ,g^{i},\varphi ^{i}\right) ,i=0,1$ satisfy $(\mathbf{H}%
_{2})$ and $V_{0}^{2}<\infty .$

$ii)$ For any $\left( t,\Pi \right) ,\varphi ^{0}\left( X\right) \leq
\varphi ^{1}\left( X\right) $ and $g^{0}\left( t,\Pi \right) \leq
g^{1}\left( t,\Pi \right) .$ Then 
\begin{equation*}
Y_{0}^{0}\leq Y_{0}^{1}.
\end{equation*}
\end{Theo}

We would like to mention that the above comparison theorem holds true only
at time $t=0$. We cannot get the result in the whole interval $\left[ 0,T%
\right] ,$ even in the Brownian case$.$ See for instance, the counterexample
which is given in \cite{WZ}.

\begin{Rema}
We should point out that\ the following cases are in fact, involved in our
present study.

\begin{enumerate}
\item \textbf{FBSDEs driven by Brownian motion}: If $\nu =0,$ then all
non--zero degree polynomials $q_{i-1}(x)$ will vanish, $H_{t}^{(1)}=W_{t}$
is a standard Brownian motion and $H_{t}^{(i)}=0,$ for $i\geq 2.$

\item \textbf{FBSDEs driven by Poisson Process}: assume that $\mu $ only has
mass at $1$, then $H_{t}^{(i)}=N_{t}-\lambda t$ is the compensated Poisson
process with intensity $\lambda $ and also $H_{t}^{(i)}=0,$ for $i\geq 2.$
For example, If we have $\nu (dx)=$\ $\sum_{j=1}^{\infty }\alpha _{j}\delta
_{\beta _{j}}(dx)$, where $\delta _{\beta _{j}}(dx)$ denotes the positive
mass measure at $\beta _{j}\in \mathbb{R}$ of size $1$. Then, The process $%
L_{\cdot }$ takes the form 
\begin{equation*}
L_{t}=at+\sum_{j=1}^{\infty }\left( N_{t}^{\left( j\right) }+\alpha
_{j}t\right) ,
\end{equation*}%
where $\{N_{t}^{\left( j\right) }\}_{j=1}^{+\infty }$ denote the sequence of
independent Poisson process with parameters $\left\{ \alpha _{j}\right\}
_{j=1}^{+\infty }.$ In this case 
\begin{equation*}
H_{t}^{(1)}=\sum_{j=1}^{\infty }\frac{\beta _{1}}{\sqrt{\alpha _{j}}}\left(
N_{t}^{\left( j\right) }+\alpha _{j}t\right)
\end{equation*}
\end{enumerate}
\end{Rema}

\section{Proofs and technical results}

\subsection{Small time duration}

In this subsection, we shall start by giving and proving the following
technical Lemma, which will be used in the proof of Theorem $\ref{STD}$. Let
us introduce the following decoupled FBSDE:%
\begin{equation}
\left\{ 
\begin{array}{ll}
\tilde{X}_{t}= & X_{0}+\int_{0}^{t}f\left( s,\tilde{X}_{s},Y_{s},Z_{s}%
\right) ds+\int_{0}^{t}\sigma \left( s,\tilde{X}_{s-},Y_{s-}\right) dH_{s},
\\ 
\tilde{Y}_{t}= & \varphi \left( X_{T}\right) +\int_{t}^{T}g\left( s,\tilde{X}%
_{s},\tilde{Y}_{s},\tilde{Z}_{s}\right) ds-\int_{t}^{T}\tilde{Z}_{s}dH_{s}.%
\end{array}%
\right.  \label{EQ2}
\end{equation}

\begin{Lemm}
\label{lem2} Assume that all the conditions in Theorem $(\ref{STD})$ are
satisfied. Let $\left( \tilde{X}_{s},\tilde{Y}_{s},\tilde{Z}_{s}\right) $
and $\left( \tilde{U}_{t},\tilde{V}_{t},\tilde{W}_{s}\right) $ belong to $%
M^{2}\left( 0,T\right) $ and satisfy the equation $(\ref{EQ2})$,\ then there
exists three constants $c,c^{\prime }$ and $c^{\prime \prime }$ depending on 
$\lambda $ and $\lambda _{0},$ such that the following estimates hold true%
\begin{equation}
\left( 1-cT^{\frac{1}{2}}\right) \mathbb{E}\sup\limits_{0\leq t\leq
T}\left\vert \tilde{X}_{t}-\tilde{U}_{t}\right\vert ^{2}\leq cT^{\frac{1}{2}%
}\left( \mathbb{E}\sup\limits_{0\leq s\leq T}\left\vert
Y_{s}-V_{s}\right\vert ^{2}+\mathbb{E}\int\limits_{0}^{T}\left\Vert
Z_{s}-W_{s}\right\Vert _{l^{2}\left( \mathbb{R}\right) }^{2}ds\right) ,
\label{A}
\end{equation}%
\begin{equation}
\left( 1-c^{\prime \prime }T\right) \mathbb{E}\left( \sup\limits_{0\leq
t\leq T}\left\vert \tilde{Y}_{t}-\tilde{V}_{t}\right\vert ^{2}\right) \leq
c^{\prime \prime }\left( 1+T\right) \mathbb{E}\left( \sup\limits_{0\leq
s\leq T}\left\vert \tilde{X}_{s}-\tilde{U}_{s}\right\vert ^{2}\right) ,
\label{B}
\end{equation}%
\begin{equation}
\mathbb{E}\left[ \int_{0}^{T}\left\Vert \tilde{Z}_{s}-\tilde{W}%
_{s}\right\Vert _{l^{2}\left( \mathbb{R}\right) }^{2}ds\right] \leq
c^{\prime }\left( \left( 1+T\right) \mathbb{E}\left( \sup\limits_{0\leq
s\leq T}\left\vert \tilde{X}_{s}-\tilde{U}_{s}\right\vert ^{2}\right) +T%
\mathbb{E}\left( \sup\limits_{0\leq s\leq T}\left\vert \tilde{Y}_{s}-\tilde{V%
}_{s}\right\vert ^{2}\right) \right) .  \label{C}
\end{equation}
\end{Lemm}

%%%%%%%%%%%%%%%%%%%%%%%%%%%%%%%%%%%%%%%%%%%%%%%%%%%%%%%%%%%
%                          Proof ofLemma <ref>lem2</ref>
%%%%%%%%%%%%%%%%%%%%%%%%%%%%%%%%%%%%%%%%%%%%%%%%%%%%%%%%%%%

\noindent \textbf{Proof of Lemma} \ref{lem2}. Let us consider $\left(
X_{t},Y_{t},Z_{t}\right) _{0\leq t\leq T},$ $\left( \tilde{X}_{t},\tilde{Y}%
_{t},\tilde{Z}_{t}\right) _{0\leq t\leq T},$

$\left( U_{t},V_{t},W_{t}\right) _{0\leq t\leq T},\left( \tilde{U}_{t},%
\tilde{V}_{t},\tilde{W}_{t}\right) _{0\leq t\leq T}\in M^{2}\left(
0,T\right) .$ First, we proceed to prove $(\ref{A}).$ Applying Itô's formula
to $\left\vert \tilde{X}_{t}-\tilde{U}_{t}\right\vert ^{2},$ taking
expectation and using the fact that $\left[ H^{i},H^{j}\right]
_{t}-\left\langle H^{i},H^{j}\right\rangle _{t}$ is an $\mathcal{F}_{t}$%
-martingale and $\left\langle H^{i},H^{j}\right\rangle _{t}=\delta _{ij}t,$
then there exists a constant $c$, depending on $\lambda $ such that%
\begin{equation*}
\begin{array}{l}
\mathbb{E}\sup\limits_{0\leq t\leq T}\left\vert \tilde{X}_{t}-\tilde{U}%
_{t}\right\vert ^{2}\leq c\left[ \mathbb{E}\int_{0}^{T}\left\vert \tilde{X}%
_{s}-\tilde{U}_{s}\right\vert \left( \left\vert \tilde{X}_{s}-\tilde{U}%
_{s}\right\vert +\left\vert Y_{s}-V_{s}\right\vert +\left\Vert
Z_{s}-W_{s}\right\Vert _{l^{2}\left( \mathbb{R}\right) }^{2}\right) ds\right.
\\ 
\left. +\mathbb{E}\int_{0}^{T}\left( \left\vert \tilde{X}_{s}-\tilde{U}%
_{s}\right\vert ^{2}+\left\vert Y_{s}-V_{s}\right\vert ^{2}\right) ds\right]
\\ 
+2\mathbb{E}\sup\limits_{0\leq t\leq T}\left\vert \int_{0}^{t}\left( \tilde{X%
}_{s-}-\tilde{U}_{s-}\right) \sigma \left( s,\tilde{X}_{s-},Y_{s-}\right)
-\sigma \left( s,\tilde{U}_{s-},V_{s-}\right) dH_{s}\right\vert .%
\end{array}%
\end{equation*}%
Burkholder-Davis-Gundy's inequality applied to the martingale%
\begin{equation*}
\int_{0}^{t}\left( \tilde{X}_{s}-\tilde{U}_{s}\right) \sigma \left( s,\tilde{%
X}_{s-},Y_{s-}\right) -\sigma \left( s,\tilde{U}_{s-},V_{s-}\right) dH_{s}
\end{equation*}%
yields the existence of a constant $C>0,$ such that%
\begin{equation*}
\begin{array}{c}
\mathbb{E}\sup\limits_{0\leq t\leq T}\left\vert \int_{0}^{t}\left( \tilde{X}%
_{s}-\tilde{U}_{s}\right) \left( \sigma \left( s,\tilde{X}%
_{s-},Y_{s-}\right) -\sigma \left( s,\tilde{U}_{s-},V_{s-}\right) \right)
dH_{s}\right\vert \\ 
\leq C\mathbb{E}\left( \left[ \int_{0}^{\cdot }\left( \tilde{X}_{s-}-\tilde{U%
}_{s-}\right) \left( \sigma \left( s,\tilde{X}_{s-},Y_{s-}\right) -\sigma
\left( s,\tilde{U}_{s-},V_{s-}\right) \right) dH_{s}\right] \right) ^{\frac{1%
}{2}}%
\end{array}%
\end{equation*}%
Moreover, since $\left\langle H^{i},H^{i}\right\rangle =\delta _{ij}t$ and $%
\left[ M\right] _{t}=\left\langle M\right\rangle _{t}+\psi _{t},$ where $%
\psi _{t}$ is a uniformly integrable martingale starting at $0$, then%
\begin{equation*}
\begin{array}{c}
\mathbb{E}\left( \sup\limits_{0\leq t\leq T}\left\vert \int_{0}^{t}\left( 
\tilde{X}_{s}-\tilde{U}_{s}\right) \left( \sigma \left( s,\tilde{X}%
_{s-},Y_{s-}\right) -\sigma \left( s,\tilde{U}_{s-},V_{s-}\right) \right)
dH_{s}\right\vert ^{2}\right) ^{1/2} \\ 
C\mathbb{E}\left( \left[ \int_{0}^{\cdot }\left( \tilde{X}_{s}-\tilde{U}%
_{s}\right) \left( \sigma \left( s,\tilde{X}_{s-},Y_{s-}\right) -\sigma
\left( s,\tilde{U}_{s-},V_{s-}\right) \right) dH_{s}\right] \right) ^{1/2}
\\ 
=C\mathbb{E}\left( \left\langle \int\limits_{0}^{\cdot }\left( \tilde{X}%
_{s-}-\tilde{U}_{s-}\right) \sigma \left( s,\tilde{X}_{s-},Y_{s-}\right)
-\sigma \left( s,\tilde{U}_{s-},V_{s-}\right) dH_{s}\right\rangle +\psi
_{t}\right) ^{1/2} \\ 
=C\mathbb{E}\left( \int\limits_{0}^{T}\left\vert \tilde{X}_{s}-\tilde{U}%
_{s}\right\vert ^{2}\left\Vert \sigma \left( s,\tilde{X}_{s},Y_{s}\right)
-\sigma \left( s,\tilde{U}_{s},V_{s}\right) \right\Vert _{l^{2}\left( 
\mathbb{R}\right) }^{2}ds\right) ^{1/2}.%
\end{array}%
\end{equation*}%
Then, modifying $c$\ if necessary, we have%
\begin{equation*}
\begin{array}{c}
\mathbb{E}\left[ \sup\limits_{0\leq t\leq T}\left\vert \tilde{X}_{t}-\tilde{U%
}_{t}\right\vert ^{2}\right] \\ 
\leq cT^{1/2}\left( \mathbb{E}\left( \sup\limits_{0\leq s\leq T}\left\vert 
\tilde{X}_{s}-\tilde{U}_{s}\right\vert ^{2}\right) +\mathbb{E}\left(
\sup\limits_{0\leq s\leq T}\left\vert Y_{s}-V_{s}\right\vert ^{2}\right) +%
\mathbb{E}\left( \int_{0}^{T}\left\Vert Z_{s}-W_{s}\right\Vert _{l^{2}\left( 
\mathbb{R}\right) }^{2}ds\right) \right) ;%
\end{array}%
\end{equation*}%
which implies that,%
\begin{equation*}
\left( 1-cT^{1/2}\right) \mathbb{E}\left( \sup\limits_{0\leq t\leq
T}\left\vert \tilde{X}_{t}-\tilde{U}_{t}\right\vert ^{2}\right) \leq
cT^{1/2}\left( \mathbb{E}\left( \sup\limits_{0\leq s\leq T}\left\vert
Y_{s}-V_{s}\right\vert ^{2}\right) +\mathbb{E}\left( \int_{0}^{T}\left\Vert
Z_{s}-W_{s}\right\Vert _{l^{2}\left( \mathbb{R}\right) }^{2}ds\right)
\right) .
\end{equation*}%
On the other hand, by applying Itô's formula to $\left\vert \tilde{Y}_{t}-%
\tilde{V}_{t}\right\vert ^{2}$, we get%
\begin{equation}
\begin{array}{l}
\left\vert \tilde{Y}_{t}-\tilde{V}_{t}\right\vert
^{2}+\int\limits_{t}^{T}\left\Vert \tilde{Z}_{s}-\tilde{W}_{s}\right\Vert
_{l^{2}\left( \mathbb{R}\right) }^{2}ds \\ 
=\left\vert \varphi \left( \tilde{X}_{T}\right) -\varphi \left( \tilde{U}%
_{T}\right) \right\vert ^{2}+2\int\limits_{t}^{T}\left( \tilde{Y}_{s}-\tilde{%
V}_{s}\right) \left( g\left( s,\tilde{X}_{s},\tilde{Y}_{s},\bar{Z}%
_{s}\right) -g\left( s,\tilde{U}_{s},\tilde{V}_{s},\tilde{W}_{s}\right)
\right) ds \\ 
-2\int\limits_{t}^{T}\left( \tilde{Y}_{s}-\tilde{V}_{s}\right) \left( \tilde{%
Z}_{s}-\tilde{W}_{s}\right) dH_{s}-\sum_{i,j}\int\limits_{t}^{T}\left( 
\tilde{Z}_{s}^{i}-\tilde{W}_{s}^{i}\right) \left( \tilde{Z}_{s}^{j}-\tilde{W}%
_{s}^{j}\right) d\left[ H^{i},H^{j}\right] _{s}.%
\end{array}
\label{ito1}
\end{equation}%
Thus, by taking expectations, invoking the assumption $\left( \mathbf{H}%
_{1}\right) $ and using the fact that $\left( \tilde{Z}_{s}^{i}-\tilde{W}%
_{s}^{i}\right) -\left\langle H^{i},H^{j}\right\rangle _{t}$ is an $\mathcal{%
F}_{t}$-martingale and $\left\langle H^{i},H^{j}\right\rangle _{t}=\delta
_{ij}t$, one can show that there exists a constant $c^{\prime }$, depending
on $\lambda $ and $\lambda _{0}$, such that%
\begin{equation*}
\begin{array}{l}
\mathbb{E}\int_{0}^{T}\left\Vert \tilde{Z}_{s}-\tilde{W}_{s}\right\Vert
_{l^{2}\left( \mathbb{R}\right) }^{2}ds\leq c^{\prime }\left[ \mathbb{E}%
\left\vert \tilde{X}_{T}-\tilde{U}_{T}\right\vert ^{2}\right. \\ 
\left. +\mathbb{E}\int_{0}^{T}\left\vert \tilde{Y}_{s}-\tilde{V}%
_{s}\right\vert \left( \left\vert \tilde{X}_{s}-\tilde{U}_{s}\right\vert
+\left\vert \tilde{Y}_{s}-\tilde{V}_{s}\right\vert +\left\Vert \tilde{Z}_{s}-%
\tilde{W}_{s}\right\Vert _{l^{2}\left( \mathbb{R}\right) }^{2}\right) ds%
\right] .%
\end{array}%
\end{equation*}%
Using the fact that $\left\vert ab\right\vert \leq \frac{1}{2}\left(
\left\vert a\right\vert ^{2}+\left\vert b\right\vert ^{2}\right) $ for any $%
a,b\in 
%TCIMACRO{\U{211d} }%
%BeginExpansion
\mathbb{R}
%EndExpansion
,$ we have%
\begin{equation*}
\begin{array}{c}
\mathbb{E}\int_{0}^{T}\left\vert \tilde{Z}_{s}-\tilde{W}_{s}\right\vert
^{2}ds\leq c^{\prime }\left[ \left( 1+T\right) \mathbb{E}\sup\limits_{0\leq
s\leq T}\left\vert \tilde{X}_{s}-\tilde{U}_{s}\right\vert ^{2}\right. \\ 
\left. +T\mathbb{E}\sup\limits_{0\leq s\leq T}\left\vert \tilde{Y}_{s}-%
\tilde{V}_{s}\right\vert ^{2}\right] +\dfrac{1}{2}\mathbb{E}%
\int_{0}^{T}\left\Vert \tilde{Z}_{s}-\tilde{W}_{s}\right\Vert _{l^{2}\left( 
\mathbb{R}\right) }^{2}ds.%
\end{array}%
\end{equation*}%
By modifying $c^{\prime }$ if necessary, we obtain%
\begin{equation}
\begin{array}{l}
\mathbb{E}\int_{0}^{T}\left\Vert \tilde{Z}_{s}-\tilde{W}_{s}\right\Vert
_{l^{2}\left( \mathbb{R}\right) }^{2}ds \\ 
\leq c^{\prime }\left[ \left( 1+T\right) \mathbb{E}\sup\limits_{0\leq s\leq
T}\left\vert \tilde{X}_{s}-\tilde{U}_{s}\right\vert ^{2}+T\mathbb{E}%
\sup\limits_{0\leq s\leq T}\left\vert \tilde{Y}_{s}-\tilde{V}_{s}\right\vert
^{2}\right] .%
\end{array}
\label{ito2}
\end{equation}%
Using equality $\left( \ref{ito1}\right) $ once again$,$ and the
Burkholder-Davis-Gundy inequality, we show that there exists a constant $%
c^{\prime \prime }$, only depending on $\lambda $ and $\lambda _{0}$, such
that%
\begin{equation*}
\begin{array}{c}
\mathbb{E}\sup\limits_{0\leq t\leq T}\left\vert \tilde{Y}_{t}-\tilde{V}%
_{t}\right\vert ^{2}\leq c^{\prime \prime }\left[ \mathbb{E}\left\vert 
\tilde{X}_{T}-\tilde{U}_{T}\right\vert ^{2}+\mathbb{E}\left(
\int_{0}^{T}\left\vert \tilde{Y}_{s}-\tilde{V}_{s}\right\vert ^{2}\left\Vert 
\tilde{Z}_{s}-\tilde{W}_{s}\right\Vert _{l^{2}\left( \mathbb{R}\right)
}^{2}ds\right) ^{1/2}\right. \\ 
\left. +\mathbb{E}\int_{0}^{T}\left\vert \tilde{Y}_{s}-\tilde{V}%
_{s}\right\vert \left( \left\vert \tilde{X}_{s}-\tilde{U}_{s}\right\vert
+\left\vert \tilde{Y}_{s}-\tilde{V}_{s}\right\vert +\left\Vert \tilde{Z}_{s}-%
\tilde{W}_{s}\right\Vert _{l^{2}\left( \mathbb{R}\right) }^{2}\right) ds%
\right] .%
\end{array}%
\end{equation*}%
Then, Taking into account $\left( \ref{ito2}\right) ,$ using Young's
inequality one more time, and modifying $c^{\prime \prime }$ if necessary,
we get 
\begin{equation*}
\begin{array}{l}
\mathbb{E}\left( \sup\limits_{0\leq t\leq T}\left\vert \tilde{Y}_{t}-\tilde{V%
}_{t}\right\vert ^{2}\right) \leq c^{\prime \prime }\left[ \left( 1+T\right) 
\mathbb{E}\left( \sup\limits_{0\leq s\leq T}\left\vert \tilde{X}_{s}-\tilde{U%
}_{s}\right\vert ^{2}\right) +T\mathbb{E}\left( \sup\limits_{0\leq s\leq
T}\left\vert \tilde{Y}_{s}-\tilde{V}_{s}\right\vert ^{2}\right) \right] \\ 
+\frac{1}{2}\mathbb{E}\left( \sup\limits_{0\leq t\leq T}\left\vert \tilde{Y}%
_{t}-\tilde{V}_{t}\right\vert ^{2}\right)%
\end{array}%
\end{equation*}%
Then, modifying $c^{\prime \prime }$ if necessary, we have%
\begin{equation*}
\left( 1-c^{^{\prime \prime }}T\right) \mathbb{E}\left( \sup\limits_{0\leq
t\leq T}\left\vert \tilde{Y}_{t}-\tilde{V}_{t}\right\vert ^{2}\right) \leq
c^{\prime \prime }\left( 1+T\right) \mathbb{E}\left( \sup\limits_{0\leq
s\leq T}\left\vert \tilde{X}_{s}-\tilde{U}_{s}\right\vert ^{2}\right) .
\end{equation*}%
Lemma \ref{lem2}\ is proved.\hfill $\square $\medskip

\noindent \textbf{Proof of Theorem }\ref{STD}. Let $\left(
X_{t},Y_{t},Z_{t}\right) _{0\leq t\leq T}$ be a possible solution of FBSDE $(%
\ref{EQ1})$ and $\left( \tilde{X},\tilde{Y},\tilde{Z}\right) $ be defined as
in Lemma $\ref{lem2}$. It is clear that the process $\tilde{X}$ is a
solution of a Forward component of the SDE $(\ref{EQ2})$, whereas the couple 
$\left( \tilde{X},\tilde{Y}\right) $ is a solution of a Backward component
of the SDE $(\ref{EQ2})$ SDE. Then $\left( \tilde{X},\tilde{Y},\tilde{Z}%
\right) $ is a solution of the above decoupled Forward Backward SDE $(\ref%
{EQ2})$. To prove the existence and the uniqueness of the solution in $%
M^{2}\left( 0,T\right) $, we use the fixed point method. Let us define a
mapping $\Psi $ from $M^{2}\left( 0,T\right) $ into itself defined by 
\begin{equation*}
\Psi \left( X,Y,Z\right) =\left( \tilde{X},\tilde{Y},\tilde{Z}\right) .
\end{equation*}%
We want to prove that there exists a constant $\delta >0,$ only depending on 
$\lambda $ and $\lambda _{0},$ such that for $T\leq \delta ,$ $\Psi $\ is a
contraction on $M^{2}\left( 0,T\right) $\ equipped with the norm%
\begin{equation*}
\begin{array}{c}
\left\Vert \Psi \left( X,Y,Z\right) \right\Vert _{M^{2}\left( 0,T\right)
}^{2}=\mathbb{E}\left\{ \sup\limits_{0\leq t\leq T}\left[ \left\vert
X_{t}\right\vert ^{2}+\left\vert Y_{t}\right\vert ^{2}\right]
+\int_{0}^{T}\left\Vert Z_{t}\right\Vert _{l^{2}\left( \mathbb{R}\right)
}^{2}dt\right\} .%
\end{array}%
\end{equation*}%
In order to achieve this goal, we firstly assume that $T\leq 1.$ Further, we
set 
\begin{equation*}
\Psi \left( X,Y,Z\right) =\left( \tilde{X},\tilde{Y},\tilde{Z}\right) ,\text{
}\Psi \left( U,V,W\right) =\left( \tilde{U},\tilde{V},\tilde{W}\right) .
\end{equation*}%
where$\left( X_{t},Y_{t},Z_{t}\right) _{0\leq t\leq T},\left(
U_{t},V_{t},W_{t}\right) _{0\leq t\leq T}$ be two elements of $M^{2}\left(
0,T\right) .$ Thus, by invoking and combining the results $\left( \ref{A}%
\right) ,$ $\left( \ref{B}\right) $ and $\left( \ref{C}\right) $ of the
Lemma $\ref{lem2},$ a simple computation shows that there exists a constant $%
\delta $ depending on $\lambda $ and $\lambda _{0},$ such that for $T\leq
\delta ,$ the following estimate holds true 
\begin{equation*}
\left\Vert \Psi \left( X,Y,Z\right) -\Psi \left( U,V,W\right) \right\Vert
_{M^{2}\left( 0,T\right) }\leq D\left\Vert \left( X,Y,Z\right) -\left(
U,V,W\right) \right\Vert _{M^{2}\left( 0,T\right) },
\end{equation*}%
For some constant $0<D<1.$This proves that the map $\Psi $ is contraction
from $M^{2}\left( 0,T\right) $ into itself$.$ Furthermore, It follows
immediately that this mapping has a unique fixed point $\left(
X_{t},Y_{t},Z_{t}\right) $ progressively measurable which is the unique
solution of FBSDE $(\ref{EQ1}).$ The proof is complete.\hfill $\square $%
\medskip

\noindent \textbf{Proof of Proposition }\ref{Prop1}$.$ Arguing as in the
proof of Lemma $\ref{lem2}$ and standard arguments of FBSDEs (see for
example \cite{A} for the Brownian case), one can prove $\left( i\right) .$
Now we proceed to prove $\left( ii\right) .$ For this end, let us define the
stopping time%
\begin{equation*}
R_{k}=\inf \left\{ t:\left\vert X_{t}\right\vert >k\right\} \text{ with }%
X_{0}=0.
\end{equation*}

For each $k,$ denoting $\left( \tilde{X},\tilde{Y},\tilde{Z}\right) =\left(
X1_{\left[ 0,R_{k}\right) },Y1_{\left[ 0,R_{k}\right) },Z1_{\left[
0,R_{k}\right) }\right) .$ It is clear that the stopped process $\tilde{X}%
=X1_{\left[ 0,R_{k}\right) }$ is bounded by $k$, and is a semimartingale as
a product of two semimartingles, which is valid for $\tilde{Y}$ as well.\
Therefore, by applying Itô's formula, using the fact that $\left[ H^{i},H^{j}%
\right] _{t}-\left\langle H^{i},H^{j}\right\rangle _{t}$ is an $\mathcal{F}%
_{t}$-martingale, $\left\langle H^{i},H^{j}\right\rangle _{t}=\delta _{ij}t$
and standart techniques from FBSDE theory, one can prove that%
\begin{equation}
\begin{array}{l}
\mathbb{E}\left\{ \sup\limits_{0\leq t\leq T}\left[ \left\vert \tilde{X}%
_{t}\right\vert ^{2p}+\left\vert \tilde{Y}_{t}\right\vert ^{2p}\right]
+\left( \int_{0}^{T}\left\Vert \tilde{Z}_{t}\right\Vert _{l^{2}\left( 
\mathbb{R}\right) }^{2}dt\right) ^{p}\right\} \\ 
\leq C_{1}\mathbb{E}\left\{ \left\vert \tilde{X}_{0}\right\vert
^{2p}+\left\vert \tilde{\varphi}\left( 0\right) \right\vert
^{2p}+\int_{0}^{T}\left[ \left\vert \tilde{f}\left( t,0,0,0\right)
\right\vert ^{2p}+\left\Vert \tilde{\sigma}\left( t,0,0\right) \right\Vert
_{l^{2}\left( \mathbb{R}\right) }^{2p}+\left\vert \tilde{g}\left(
t,0,0,0\right) \right\vert ^{2p}\right] dt\right\} \\ 
+\dsum\limits_{0<s\leq t}\left\{ \left( \tilde{X}_{s}^{2p}\right) -\tilde{X}%
_{s-}^{2p}-2p\tilde{X}_{s}^{2p-1}\Delta \tilde{X}_{s}-p\left( 2p-1\right) 
\text{ }\tilde{X}_{s}^{2p-2}\left( \Delta \tilde{X}_{s}\right) ^{2}\right\}
\\ 
+\dsum\limits_{t<s\leq T}\left\{ \tilde{Y}_{s}^{2p}-\tilde{Y}_{s-}^{2p}-2p%
\tilde{Y}_{s}^{2p-1}\Delta \tilde{Y}_{s}-p\left( 2p-1\right)
Y_{s}^{2p-2}\left( \Delta \tilde{Y}_{s}\right) ^{2}\right\} ,%
\end{array}
\label{x4}
\end{equation}

where we have denoted by $\tilde{\varphi},$ $\tilde{f},$ $\tilde{\sigma},$
and $\tilde{g}$ the restriction of the functions of $\varphi ,$ $f,$ $\sigma
,$ and $g.$ Now, we proceed to prove that 
\begin{equation*}
\dsum\limits_{0<s\leq t}\left\{ X_{s}^{2p}-X_{s-}^{2p}-2pX_{s}^{2p-1}\Delta
X_{s}-p\left( 2p-1\right) X_{s}^{2p-2}\left( \Delta X_{s}\right)
^{2}\right\} <C\left[ X,X\right] _{t}.
\end{equation*}

Since $\tilde{X}$ takes its values in intervals of the form $\left[ -k,k%
\right] $, for $h\left( x\right) =$ $x^{2p},$ it is easy to show that 
\begin{equation*}
\left\vert h\left( x\right) -h\left( y\right) -\left( y-x\right) h^{^{\prime
}}\left( x\right) -\left( y-x\right) ^{2}h^{^{\prime \prime }}\left(
x\right) \right\vert \leq C\left( y-x\right) ^{2}
\end{equation*}

Thus%
\begin{equation}
\begin{array}{l}
\dsum\limits_{0<s\leq t}\left\vert \tilde{X}_{s}^{2p}-\tilde{X}_{s-}^{2p}-2p%
\tilde{X}_{s}^{2p-1}\Delta \tilde{X}_{s}-p\left( 2p-1\right) \tilde{X}%
_{s}^{2p-2}\left( \Delta \tilde{X}_{s}\right) ^{2}\right\vert \\ 
\leq C\dsum\limits_{0<s\leq t}^{2}\left( \Delta \tilde{X}_{s}\right) \leq C%
\left[ \tilde{X},\tilde{X}\right] _{t}<\infty .%
\end{array}
\label{sum1}
\end{equation}%
Therefore by similar arguments developed above, one can easily derive that%
\begin{equation}
\begin{array}{l}
\dsum\limits_{t<s\leq T}\left\vert \tilde{Y}_{s}^{2p}-\tilde{Y}_{s-}^{2p}-2p%
\tilde{Y}_{s}^{2p-1}\Delta \tilde{Y}_{s}-p\left( 2p-1\right) \tilde{Y}%
_{s}^{2p-2}\left( \Delta \tilde{Y}_{s}\right) ^{2}\right\vert \\ 
\leq C\dsum\limits_{t<s\leq T}\left( \Delta \tilde{Y}_{s}\right) ^{2}\leq
C\left( \left[ \tilde{Y},\tilde{Y}\right] _{T}-\left[ \tilde{Y},\tilde{Y}%
\right] _{t}\right) <\infty .%
\end{array}
\label{sum2}
\end{equation}

Combining $(\ref{sum1}),$ $(\ref{sum2})$ and $(\ref{x4}),$ we get%
\begin{equation*}
\begin{array}{l}
\mathbb{E}\left\{ \sup\limits_{0\leq t\leq T}\left[ \left\vert \tilde{X}%
_{t}\right\vert ^{2p}+\left\vert \tilde{Y}_{t}\right\vert ^{2p}\right]
+\left( \int_{0}^{T}\left\Vert \tilde{Z}_{t}\right\Vert _{l^{2}\left( 
\mathbb{R}\right) }^{2}dt\right) ^{p}\right\} \\ 
\leq C_{1}\mathbb{E}\left\{ \left\vert \tilde{X}_{0}\right\vert
^{2p}+\left\vert \varphi \left( 0\right) \right\vert ^{2p}+\int_{0}^{T}\left[
\left\vert f\left( t,0,0,0\right) \right\vert ^{2p}+\left\Vert \tilde{\sigma}%
\left( t,0,0\right) \right\Vert _{l^{2}\left( \mathbb{R}\right)
}^{2p}+\left\vert g\left( t,0,0,0\right) \right\vert ^{2p}\right] dt\right\}
+C<\infty .%
\end{array}%
\end{equation*}%
Since the last inequality is valid for $\left( \tilde{X},\tilde{Y},\tilde{Z}%
\right) $ for each $k$, it also remains valid for $\left( X,Y,Z\right) $ and
this completes the proof.\hfill $\square $

\subsection{Large time duration}

To prove Theorem $\ref{LTD}$, we need the following proposition, which
allows us to prove global existence and uniqueness of the equation $(\ref%
{EQ1})$. By using similar arguments introduced in \cite{Z} consisting in
solving the system iteratively in small intervals having fixed length.

\begin{Prop}
\label{Prop2} Let $\Pi ^{i},i=0,1,$ be the solution to FBSDEs:%
\begin{equation*}
\left\{ 
\begin{array}{ll}
X_{t}^{i}= & x_{i}+\int_{0}^{t}f\left( s,\Pi _{s}^{i}\right)
ds+\int_{0}^{t}\sigma \left( s,X_{s-}^{i},Y_{s-}^{i}\right) dH_{s}, \\ 
Y_{t}^{i}= & \varphi \left( X_{T}^{i}\right) +\int_{t}^{T}g\left( s,\Pi
_{s}^{i}\right) ds-\int_{t}^{T}Z_{s}^{i}dH_{s}.%
\end{array}%
\right.
\end{equation*}%
Assume that\ $\left( \mathbf{H}_{1}\right) $ is satisfied and $%
V_{0}^{2}<\infty .$ Then%
\begin{equation*}
\left\vert Y_{0}^{1}-Y_{0}^{0}\right\vert \leq \bar{\lambda}_{0}\left\vert
x_{1}-x_{0}\right\vert ,
\end{equation*}%
where%
\begin{equation}
\bar{\lambda}_{0}=c\left( \left[ \lambda _{0}+1\right] e^{\left( 2\lambda
+\lambda ^{2}\right) T}-1\right) .  \label{Lamda}
\end{equation}
\end{Prop}

The following lemma gives estimates of $\bar{\lambda}_{0}$ in terms of $%
\lambda $ and $\lambda _{0}$. This estimation is the key step for the proof
of Theorem \ref{LTD}$.$

\begin{Lemm}
\label{lem3} Consider the following linear FBSDE:%
\begin{equation}
\left\{ 
\begin{array}{ll}
X_{t}= & 1+\int_{0}^{t}\left(
a_{s}^{1}X_{s}+b_{s}^{1}Y_{s}+c_{s}^{1}Z_{s}\right) ds+\int_{0}^{t}\left(
a_{s}^{2}X_{s}+b_{s}^{2}Y_{s}\right) dH_{s}, \\ 
Y_{t}= & FX_{T}+\int_{t}^{T}\left(
a_{s}^{3}X_{s}+b_{s}^{3}Y_{s}+c_{s}^{3}Z_{s}\right)
ds-\int_{t}^{T}Z_{s}dH_{s}.%
\end{array}%
\right.  \label{EQ33}
\end{equation}%
Assume $\left\vert a_{t}^{i}\right\vert ,\left\vert b_{t}^{i}\right\vert
,\left\vert c_{t}^{i}\right\vert \leq \lambda ,i=1,2,3$ and $\left\vert
F\right\vert \leq \lambda _{0}.$ Let $\delta $ be as in theorem \ref{STD}$.$
And assume further that 
\begin{equation}
b_{t}^{2}c_{t}^{1}=0;\text{ \ \ }%
b_{t}^{1}+a_{t}^{2}c_{t}^{1}+b_{t}^{2}c_{t}^{3}=0.  \label{H.3}
\end{equation}%
Then for $T\leq \delta ,$

$i)$ The LFBSDE $(\ref{EQ33})$ admits a unique solution.

$ii)$ 
\begin{equation}
\left\vert Y_{0}\right\vert \leq \bar{\lambda}_{0},  \label{y0}
\end{equation}%
where $\bar{\lambda}_{0}$\ is defined by $(\ref{Lamda}).$
\end{Lemm}

\noindent \textbf{Proof of Lemma} \ref{lem3}$.$ First, we can easily check
that LFBSDE $(\ref{EQ33})$ satisfy assumptions of Theorem $\ref{STD}$, then
it has a unique solution $\left( X_{t},Y_{t},Z_{t}\right) $ which belongs to
the space $M^{2}\left( 0,T\right) $. This gives the proof of the assertion $%
\left( i\right) .$

We shall prove the assertion $\left( ii\right) $. We split the proof into
two steps.\medskip

$Step1$. For any $t\in $\ $\left[ 0,T\right) $ and any $\xi \in L^{2}\left( 
\mathcal{F}_{0}\right) ,$ we put $\bar{\Pi}_{s}\overset{\bigtriangleup }{=}%
\left( X_{t}\xi ,Y_{t}\xi ,Z_{t}\xi \right) ,$ $s\in \left[ t,T\right] .$
Then $\bar{\Pi}_{s}$\ satisfies the following linear FBSDE 
\begin{equation*}
\left\{ 
\begin{array}{ll}
\bar{X}_{s}= & X_{t}\xi +\int\limits_{t}^{s}\left[ a_{r}^{1}\bar{X}%
_{r}+b_{r}^{1}\bar{Y}_{r}+c_{r}^{1}\bar{Z}_{r}\right] dr+\int\limits_{t}^{s}%
\left[ a_{r}^{2}\bar{X}_{r}+b_{r}^{2}\bar{Y}_{r}\right] dH_{r}, \\ 
\bar{Y}_{s}= & F\bar{X}_{T}+\int\limits_{s}^{T}\left[ a_{r}^{3}\bar{X}%
_{r}+b_{r}^{3}\bar{Y}_{r}+c_{r}^{3}\bar{Z}_{r}\right] dr-\int\limits_{s}^{T}%
\bar{Z}_{r}dH_{r}.%
\end{array}%
\right.
\end{equation*}%
By assertion $(ii)$ of Theorem $\ref{STD}$, we get%
\begin{equation*}
E\left\{ \left\vert \bar{Y}_{t}\right\vert ^{2}\right\} =E\left\{ \left\vert
Y_{t}\xi \right\vert ^{2}\right\} \leq C_{0}^{2}E\left\{ \left\vert X_{t}\xi
\right\vert ^{2}\right\} .
\end{equation*}%
Since $\xi $ is arbitrary, we have $\left\vert Y_{t}\right\vert \leq $\ $%
C_{0}\left\vert X_{t}\right\vert $\ , $P$-a.s.,$\forall t.$\medskip

$Step2$. We define%
\begin{equation*}
\tau \overset{\Delta }{=}\inf \left\{ t>0:X_{t}=0\right\} \wedge T;\text{ \
and\ }\tau _{n}\overset{\Delta }{=}\inf \left\{ t>0:X_{t}=\frac{1}{n}%
\right\} \wedge T.
\end{equation*}%
Then $\tau _{n}\uparrow \tau $ and $X_{t}>0$ for $t\in $\ $\left[ 0,\tau
\right) .$ We also define the pure jump process $\eta ,$ by the following
formula%
\begin{equation*}
\eta _{t}=\dprod\limits_{0<s\leq t}\left( 1-\left( X_{s}\right) ^{-1}\Delta
X_{s}\right) \frac{\left( X_{s-}\right) ^{-1}}{\left( X_{s}\right) ^{-1}}
\end{equation*}%
The above product is clearly càdlàg, adapted, converges and is of finite
variation. We put for any $t\in $\ $\left[ 0,\tau \right) ,$%
\begin{equation*}
A_{t}=\eta _{t}\left( X_{t}\right) ^{-1}.
\end{equation*}%
It should be noted that when we apply Itô's formula to $\left( X_{t}\right)
^{-1},$ a sum of discontinuous quantities appears. To eliminate this, we
shall apply Itô's formula to $A_{t}=\eta _{t}\left( X_{t}\right) ^{-1}$
instead of $\left( X_{t}\right) ^{-1}.$ Firstly, applying Itô's formula to $%
A_{t}$, we have%
\begin{equation}
\begin{array}{c}
A_{t}=A_{0}-\int_{0}^{t}\eta _{s-}\left( X_{s-}\right)
^{-2}dX_{s}+\int_{0}^{t}\left( X_{s-}\right) ^{-1}d\eta
_{s}+\int_{0}^{t}A_{s-}\left( X_{s}\right) ^{-2}d\left[ X,X\right] _{s}^{c}
\\ 
+\dsum\limits_{0<s\leq t}\left( A_{s}-A_{s-}+A_{s-}\left( X_{s-}\right)
^{-1}\left( \Delta X_{s}\right) -\left( X_{s-}\right) ^{-1}\Delta \eta
_{s}\right) ,%
\end{array}
\label{At}
\end{equation}%
Note that $\eta $ is a pure jump process. Hence $\left[ \eta ,X\right] ^{c}=%
\left[ \eta ,\eta \right] ^{c}=0$ and

\begin{equation*}
\int_{0}^{t}\left( \tilde{X}_{s-}\right) ^{-1}d\eta _{s}=\sum_{0<s\leq
t}\left( \tilde{X}_{s-}\right) ^{-1}\Delta \eta _{s}.
\end{equation*}%
Then $(\ref{At})$ becomes%
\begin{equation*}
\begin{array}{c}
A_{t}=A_{0}-\int_{0}^{t}\eta _{s-}\left( X_{s-}\right)
^{-2}dX_{s}+\int_{0}^{t}A_{s-}\left( X_{s-}\right) ^{-2}d\left[ X,X\right]
_{s}^{c} \\ 
+\dsum\limits_{0<s\leq t}\left( A_{s}-A_{s-}+A_{s-}\left( X_{s-}\right)
^{-1}\Delta X_{s}\right) ,%
\end{array}%
\end{equation*}%
The following equality is obvious, from the definition of the process $A,$ 
\begin{equation*}
A_{s}=A_{s-}\left( 1-\left( X_{t}\right) ^{-1}\Delta X_{t}\right) .
\end{equation*}%
Now by replacing the above equality into the previous one, one can get%
\begin{equation*}
\sum_{0<s\leq t}\left( A_{s}-A_{s-}+A_{s-}\left( X_{s-}\right) ^{-1}\Delta
X_{s}\right) =0.
\end{equation*}%
Therefore,%
\begin{equation*}
A_{t}=A_{0}-\int_{0}^{t}A_{s}\left( X_{s}\right)
^{-1}dX_{s}+\int_{0}^{t}A_{s-}\left( X_{s}\right) ^{-2}d\left[ X,X\right]
_{s}^{c},
\end{equation*}%
with%
\begin{equation*}
d\left[ X,X\right] _{s}^{c}=\sum_{i,j}\left(
a_{s}^{2,i}X_{s}+b_{s}^{2,i}Y_{s}\right) \left(
a_{s}^{2,j}X_{s}+b_{s}^{2,j}Y_{s}\right) q_{i-1}\left( 0\right)
q_{j-1}\left( 0\right) ds.
\end{equation*}%
Thanks to Lemma $\ref{lem}$, we get%
\begin{eqnarray*}
d\left[ X,X\right] _{s}^{c} &=&\left[ \left(
a_{s}^{2}X_{s}+b_{s}^{2}Y_{s}\right) ^{2}-\sum_{i,j}\left(
a_{s}^{2,i}X_{s}+b_{s}^{2,i}Y_{s}\right) \left(
a_{s}^{2,j}X_{s}+b_{s}^{2,j}Y_{s}\right) \int_{%
%TCIMACRO{\U{211d} }%
%BeginExpansion
\mathbb{R}
%EndExpansion
}p_{i}\left( x\right) p_{j}\left( x\right) v\left( dx\right) \right] ds. \\
&=&\left[ \left( a_{s}^{2}X_{s}+b_{s}^{2}Y_{s}\right) ^{2}-\Psi _{s}\right]
ds
\end{eqnarray*}%
Hence%
\begin{equation*}
\begin{array}{l}
A_{t}=A_{0}-\int_{0}^{t}\left[ A_{s}\left( X_{s}\right) ^{-1}\left(
a_{s}^{1}X_{s}+b_{s}^{1}Y_{s}+c_{s}^{1}Z_{s}\right) -A_{s}\left(
X_{s}\right) ^{-2}\left( a_{s}^{2}X_{s}+b_{s}^{2}Y_{s}\right) ^{2}\right] ds
\\ 
-\int_{0}^{t}A_{s-}\left( X_{s-}\right) ^{-1}\left(
a_{s}^{2}X_{s-}+b_{s}^{2}Y_{s-}\right) dH_{s}-\int_{0}^{t}A_{s}\left(
X_{s}\right) ^{-2}\Psi _{s}ds.%
\end{array}%
\end{equation*}%
Let us define the following processes%
\begin{equation*}
\hat{Y}_{t}=Y_{t}A_{t}\text{; \ \ }\hat{Z}_{t}\overset{\bigtriangleup }{=}%
A_{t}Z_{t}-A_{t}\left( X_{t}\right) ^{-1}Y_{t}\left(
a_{t}^{2}X_{t}+b_{t}^{2}Y_{t}\right) .
\end{equation*}%
Then after the result of the\ Step $1$, we have%
\begin{equation*}
\left\vert \hat{Y}_{t}\right\vert \leq C_{0}.
\end{equation*}%
Now, applying Itô's formula to $\hat{Y}_{t}$, we obtain%
\begin{equation*}
\begin{array}{l}
d\hat{Y}_{t}=-A_{t}\left(
a_{t}^{3}X_{t}+b_{t}^{3}Y_{t}+c_{t}^{3}Z_{t}\right) dt \\ 
-\left[ Y_{t}A_{t}\left( X_{t}\right) ^{-1}\left(
a_{t}^{1}X_{t}+b_{t}^{1}Y_{t}+c_{t}^{1}Z_{t}\right) dt-A_{t}\left(
X_{t}\right) ^{-2}\left( a_{t}^{2}X_{t}+b_{t}^{2}Y_{t}\right) ^{2}\right] dt
\\ 
-\left[ A_{t}\left( X_{t}\right) ^{-1}\left(
a_{t}^{2}X_{t}+b_{t}^{2}Y_{t}\right) Z_{t}\right] dt-\left[ Y_{t}A_{t}\left(
X_{t}\right) ^{-2}\Psi _{t}\right] dt \\ 
+\left[ A_{t-}Z_{t}-Y_{t-}A_{t-}\left( X_{t-}\right) ^{-1}\left(
a_{t}^{2}X_{t-}+b_{t}^{2}Y_{t-}\right) \right] dH_{t}+d\tilde{A}_{t}.%
\end{array}%
\end{equation*}%
where we have denoted by $\tilde{A}_{t}=\left[ A,Y\right] _{t}-\left\langle
A,Y\right\rangle _{t}.$ By using the definition of the processes$\left( \hat{%
Y},\hat{Z}\right) $ it follows that%
\begin{equation*}
\begin{array}{l}
d\hat{Y}_{t}=\hat{Z}_{t}dH_{t}-\left[ c_{t}^{3}+c_{t}^{1}\eta _{t}^{-1}\hat{Y%
}_{t}+a_{t}^{2}+b_{t}^{2}\eta _{t}^{-1}\hat{Y}_{t}\right] \hat{Z}_{t}dt \\ 
-\left[ c_{t}^{1}b_{t}^{2}\left( \eta _{t}^{-1}\right) ^{2}\hat{Y}%
_{t}^{3}+\left( b_{t}^{1}+a_{t}^{2}c_{t}^{1}+c_{t}^{3}b_{t}^{2}\right) \eta
_{t}^{-1}\hat{Y}_{t}^{2}\right] dt-\left[ a_{t}^{3}\eta _{t}+\left(
b_{t}^{3}+a_{t}^{1}+c_{t}^{3}a_{t}^{2}\right) \hat{Y}_{t}\right] dt \\ 
-\left[ Y_{t}A_{t}\left( X_{t}\right) ^{-2}\Psi _{t}\right] dt+d\tilde{A}%
_{t}.%
\end{array}%
\end{equation*}%
Thus, by taking into account $(\ref{H.3}),$ 
\begin{equation*}
\begin{array}{l}
d\hat{Y}_{t}=\hat{Z}_{t}dH_{t}-\left[ c_{t}^{3}+c_{t}^{1}\eta _{t}^{-1}\hat{Y%
}_{t}+a_{t}^{2}+b_{t}^{2}\eta _{t}^{-1}\hat{Y}_{t}\right] \hat{Z}_{t}dt \\ 
-\left[ a_{t}^{3}\eta _{t}+\left(
b_{t}^{3}+a_{t}^{1}+c_{t}^{3}a_{t}^{2}\right) \hat{Y}_{t}\right] dt-\left[
Y_{t}A_{t}\left( X_{t}\right) ^{-2}\Psi _{t}\right] dt+d\tilde{A}_{t}.%
\end{array}%
\end{equation*}%
We put%
\begin{equation*}
\Gamma _{t}=1+\int_{0}^{t}\Gamma _{s}\left( X_{s}\right) ^{-2}\Psi
_{s}1_{\left\{ \tau >s\right\} }ds.
\end{equation*}%
\begin{equation*}
M_{t}=1+\dsum\limits_{i=1}^{\infty }\int_{0}^{t}\left( q_{i-1}\left(
0\right) \right) ^{-1}M_{s}\left( \left( c_{s}^{3}+a_{s}^{2}\right) +\left(
c_{s}^{1}+b_{s}^{2}\right) \eta _{s}^{-1}\hat{Y}_{s}\right) 1_{\left\{ \tau
>s\right\} }dB_{s};
\end{equation*}%
\begin{equation*}
N_{t}=1+\int_{0}^{t}N_{s}\left(
a_{s}^{1}+b_{s}^{3}+a_{s}^{2}c_{s}^{3}\right) 1_{\left\{ \tau >s\right\} }ds.
\end{equation*}%
Applying Itô's formula to $\left( \Gamma _{t}N_{t}M_{t}\hat{Y}_{t}\right) ,$
we obtain%
\begin{equation*}
\begin{array}{l}
d\left( \Gamma _{t}N_{t}M_{t}\hat{Y}_{t}\right) =\Gamma _{t}N_{t}M_{t}\hat{Z}%
_{t}1_{\left\{ \tau >t\right\} }dH_{t}+\dsum\limits_{i=1}^{\infty }\left(
\Gamma _{t}\hat{Y}_{t}N_{t}M_{t}\left( q_{i-1}\left( 0\right) \right)
^{-1}\left( c_{t}^{3}+a_{t}^{2}\right) +\left( c_{t}^{1}+b_{t}^{2}\right)
\eta _{t}^{-1}\hat{Y}_{t}\right) 1_{\left\{ \tau >t\right\} }dB_{t} \\ 
-\eta _{t}\Gamma _{t}N_{t}M_{t}a_{t}^{3}1_{\left\{ \tau >s\right\}
}dt+\Gamma _{t}M_{t}N_{t}d\tilde{A}_{t}.%
\end{array}%
\end{equation*}%
Taking expectations, we get%
\begin{equation}
Y_{0}=E\left( \Gamma _{\tau _{n}}N_{\tau _{n}}M_{\tau _{n}}\hat{Y}_{\tau
_{n}}+\int_{0}^{\tau _{n}}\eta _{t}\Gamma _{t}N_{t}M_{t}a_{t}^{3}dt\right) .
\label{Y0}
\end{equation}%
Since $\left\vert \hat{Y}_{t}\right\vert \leq C_{0},$ $M$ is an $\mathcal{F}%
_{t}$--martingale and $\left\vert N_{t}\right\vert \leq e^{\left( 2\lambda
+\lambda ^{2}\right) t}.$ Moreover, we observe that

if $\tau =T,$ $\left\vert Y_{\tau }\right\vert =\left\vert Y_{T}\right\vert
=\left\vert FX_{T}\right\vert =\left\vert FX_{\tau }\right\vert \leq \lambda
_{0}\left\vert X_{\tau }\right\vert .$

If $\tau <T,$ $X_{\tau }=0,$ and thus $\left\vert Y_{\tau }\right\vert \leq
C_{0}\left\vert X_{\tau }\right\vert =0.$

Therefore, in both cases it holds that $\left\vert Y_{\tau }\right\vert \leq
\lambda _{0}\left\vert X_{\tau }\right\vert .$

\noindent Now, applying Ito's formula to $\left\vert Y_{t}\right\vert ^{2}$
from $s=\tau _{n}$ to $s=\tau ,$ we obtain%
\begin{equation*}
\begin{array}{l}
\left\vert Y_{\tau _{n}}\right\vert ^{2}+E_{\tau _{n}}\left( \int_{\tau
_{n}}^{\tau }\left\Vert Z_{t}\right\Vert _{l^{2}\left( \mathbb{R}\right)
}^{2}dt\right) =E_{\tau _{n}}\left( \left\vert Y_{\tau }\right\vert
^{2}+2\int_{\tau _{n}}^{\tau }Y_{t}\left(
a_{t}^{3}X_{t}+b_{t}^{3}Y_{t}+c_{t}^{3}Z_{t}\right) dt\right) \\ 
\leq E_{\tau _{n}}\left( \lambda _{0}^{2}\left\vert X_{\tau }\right\vert
^{2}+C\int_{\tau _{n}}^{\tau }\left( \left\vert X_{t}\right\vert
^{2}+\left\vert Y_{t}\right\vert ^{2}\right) dt+\dfrac{1}{2}\int_{\tau
_{n}}^{\tau }\left\Vert Z_{t}\right\Vert _{l^{2}\left( \mathbb{R}\right)
}^{2}dt\right) .%
\end{array}%
\end{equation*}%
Similarly, applying Ito's formula to $\left\vert X_{t}\right\vert ^{2}$ from 
$s=\tau _{n}$ to $s=\tau ,$ we obtain,%
\begin{equation*}
E_{\tau _{n}}\left( \left\vert X_{\tau }\right\vert ^{2}\right) \leq E_{\tau
_{n}}\left( \left\vert X_{\tau _{n}}\right\vert ^{2}+C\int_{\tau _{n}}^{\tau
}\left( \left\vert X_{t}\right\vert ^{2}+\left\vert Y_{t}\right\vert
^{2}\right) dt+\frac{1}{2\lambda _{0}^{2}}\int_{\tau _{n}}^{\tau }\left\Vert
Z_{t}\right\Vert _{l^{2}\left( \mathbb{R}\right) }^{2}dt\right) .
\end{equation*}%
Thus%
\begin{equation*}
\left\vert Y_{\tau _{n}}\right\vert ^{2}\leq E_{\tau _{n}}\left( \lambda
_{0}^{2}\left\vert X_{\tau _{n}}\right\vert ^{2}+C\int_{\tau _{n}}^{\tau
}\left( \left\vert X_{t}\right\vert ^{2}+\left\vert Y_{t}\right\vert
^{2}\right) dt\right) .
\end{equation*}%
Note that $\left\vert X_{\tau _{n}}\right\vert \geq \frac{1}{n},$ then%
\begin{equation*}
\begin{array}{l}
\left\vert \hat{Y}_{\tau _{n}}\right\vert \leq \lambda _{0}\left\vert \eta
_{\tau _{n}}\right\vert +CE_{\tau _{n}}^{\frac{1}{2}}\left( \int_{\tau
_{n}}^{\tau }\left( \left\vert \tilde{X}_{t}\right\vert ^{2}+\left\vert 
\tilde{Y}_{t}\right\vert ^{2}\right) dt\right) \\ 
\leq \lambda _{0}\left\vert \eta _{\tau _{n}}\right\vert +CE_{\tau _{n}}^{%
\frac{1}{2}}\left( \sup\limits_{\tau _{n}\leq t\leq \tau }\left( \left\vert 
\tilde{X}_{t}\right\vert ^{2}+\left\vert \tilde{Y}_{t}\right\vert
^{2}\right) \left( \tau -\tau _{n}\right) \right) ,%
\end{array}%
\end{equation*}%
where%
\begin{equation*}
\tilde{X}_{t}\overset{\bigtriangleup }{=}X_{t}\left\vert \eta _{\tau
_{n}}\right\vert \left( X_{\tau _{n}}\right) ^{-1}\ ;\text{ \ }\tilde{Y}_{t}%
\overset{\bigtriangleup }{=}Y_{t}\left\vert \eta _{\tau _{n}}\right\vert
\left( X_{\tau _{n}}\right) ^{-1}.
\end{equation*}%
Now by $(\ref{Y0})$, we get%
\begin{equation*}
\begin{array}{l}
\left\vert \hat{Y}_{0}\right\vert \leq \lambda E\left( \Gamma
_{t}M_{t}\right) \int_{0}^{T}\left\vert \eta _{t}\right\vert e^{\left(
2\lambda +\lambda ^{2}\right) t}dt \\ 
+E\left\{ e^{\left( 2\lambda +\lambda ^{2}\right) T}M_{\tau _{n}}\Gamma
_{\tau _{n}}\left( \left\vert \eta _{\tau _{n}}\right\vert \lambda
_{0}+CE_{\tau _{n}}^{\frac{1}{2}}\left( \sup\limits_{\tau _{n}\leq t\leq
\tau }\left( \left\vert \tilde{X}_{t}\right\vert ^{2}+\left\vert \tilde{Y}%
_{t}\right\vert ^{2}\right) \left( \tau -\tau _{n}\right) \right) \right)
\right\} \\ 
\leq c^{\prime }\left( e^{\left( 2\lambda +\lambda ^{2}\right) T}-1\right)
+c^{\prime \prime }\lambda _{0}e^{\left( 2\lambda +\lambda ^{2}\right) T} \\ 
+CE\left\{ M_{\tau _{n}}\Gamma _{\tau _{n}}E_{\tau _{n}}^{\frac{1}{2}}\left(
\sup\limits_{\tau _{n}\leq t\leq \tau }\left( \left\vert \tilde{X}%
_{t}\right\vert ^{2}+\left\vert \tilde{Y}_{t}\right\vert ^{2}\right) \left(
\tau -\tau _{n}\right) \right) \right\} \\ 
\leq \bar{\lambda}_{0}+CE^{\frac{1}{2}}\left( \left\vert M_{\tau
_{n}}\right\vert ^{2}\left\vert \Gamma _{\tau _{n}}\right\vert ^{2}\right)
E^{\frac{1}{2}}\left( \sup\limits_{\tau _{n}\leq t\leq \tau }\left(
\left\vert \tilde{X}_{t}\right\vert ^{2}+\left\vert \tilde{Y}_{t}\right\vert
^{2}\right) \left( \tau -\tau _{n}\right) \right) \\ 
\leq \bar{\lambda}_{0}+CE^{\frac{1}{4}}\left( \sup\limits_{\tau _{n}\leq
t\leq \tau }\left( \left\vert \tilde{X}_{t}\right\vert ^{4}+\left\vert 
\tilde{Y}_{t}\right\vert ^{4}\right) \right) E^{\frac{1}{4}}\left(
\left\vert \tau -\tau _{n}\right\vert ^{2}\right) .%
\end{array}%
\end{equation*}%
Note that $\left( \tilde{X}_{t},\tilde{Y}_{t}\right) $ satisfies the
following LFBSDE:%
\begin{equation*}
\left\{ 
\begin{array}{ll}
\tilde{X}_{t}= & 1+\int_{0}^{t}\left[ a_{r}^{1}1_{\left\{ \tau _{n}\leq
r\right\} }\tilde{X}_{r}+b_{r}^{1}1_{\left\{ \tau _{n}\leq r\right\} }\tilde{%
Y}_{r}+c_{r}^{1}1_{\left\{ \tau _{n}\leq r\right\} }\tilde{Z}_{r}\right] dr
\\ 
& +\int_{0}^{t}\left[ a_{r}^{2}1_{\left\{ \tau _{n}\leq r\right\} }\tilde{X}%
_{r}+b_{r}^{2}1_{\left\{ \tau _{n}\leq r\right\} }\tilde{Y}_{r}\right]
dH_{r}, \\ 
\tilde{Y}_{t}= & F\tilde{X}_{T}+\int_{t}^{T}\left[ a_{r}^{3}1_{\left\{ \tau
_{n}\leq r\right\} }\tilde{X}_{r}+b_{r}^{3}1_{\left\{ \tau _{n}\leq
r\right\} }\tilde{Y}_{r}+c_{r}^{3}1_{\left\{ \tau _{n}\leq r\right\} }\tilde{%
Z}_{r}\right] dr-\int_{t}^{T}\tilde{Z}_{r}dH_{r}.%
\end{array}%
\right.
\end{equation*}%
By $(ii)$ of Proposition $\ref{Prop1}$, we have%
\begin{equation*}
E\left( \sup\limits_{\tau _{n}\leq t\leq \tau }\left( \left\vert \tilde{X}%
_{t}\right\vert ^{4}+\left\vert \tilde{Y}_{t}\right\vert ^{4}\right) \right)
\leq E\left( \sup\limits_{0\leq t\leq T}\left( \left\vert \tilde{X}%
_{t}\right\vert ^{4}+\left\vert \tilde{Y}_{t}\right\vert ^{4}\right) \right)
\leq C_{1}.
\end{equation*}%
Thus%
\begin{equation*}
\left\vert \hat{Y}_{0}\right\vert \leq \bar{\lambda}_{0}+CE^{\frac{1}{4}%
}\left( \left\vert \tau -\tau _{n}\right\vert ^{2}\right) .
\end{equation*}%
Then for $n\rightarrow \infty ,$ we get $\left\vert \hat{Y}_{0}\right\vert
\leq \bar{\lambda}_{0}.$That is, $\left\vert Y_{0}\right\vert \leq \bar{%
\lambda}_{0}\left\vert X_{0}\right\vert \left\vert \eta _{0}\right\vert =%
\bar{\lambda}_{0}.$This complete the proof.\hfill $\square $\bigskip

\textbf{Proof of Proposition }$\ref{Prop2}.$ The proof is the same as in
Corollary 1 in $\cite{Z}$, by replacing the Brownian part by the Teugels
martingales and using the above lemma.\hfill $\square $\bigskip

Now we are able to give the proof of our main result. We shall extend by
induction the theorem $\ref{STD}$ to $\ref{LTD}$.\bigskip

\textbf{Proof of Theorem }$\ref{LTD}.$ First we prove $\left( i\right) $.
Let $\lambda $\ and $\lambda _{0}$\ be as in Theorem $\ref{STD},$\ and $\bar{%
\lambda}_{0}$\ is a constant defined as in $(\ref{Lamda})$. Let $\delta $\
be a constant as in Theorem $\ref{STD},$\ but corresponding to $\lambda $\
and $\bar{\lambda}_{0}$\ instead of $\lambda $\ and $\lambda _{0}$. For some
integer $n$, we assume\textbf{\ }$\left( n-1\right) \delta <T\leq n\delta $%
\textbf{\ }and consider a partition of $\left[ 0,T\right] ,$ with $T_{i}%
\overset{\vartriangle }{=}\frac{iT}{n},i=0,...,n$.

We consider the mapping:%
\begin{equation*}
\begin{array}{ccc}
G_{n}: & \Omega \times 
%TCIMACRO{\U{211d} }%
%BeginExpansion
\mathbb{R}
%EndExpansion
\rightarrow & 
%TCIMACRO{\U{211d} }%
%BeginExpansion
\mathbb{R}
%EndExpansion
\\ 
& \omega \times x\mapsto & \varphi \left( \omega ,x\right)%
\end{array}%
\end{equation*}%
Let us consider the following FBSDE over the small interval $\left[
T_{n-1},T_{n}\right] $,%
\begin{equation}
\left\{ 
\begin{array}{ll}
X_{t}^{n}= & x+\int_{T_{n-1}}^{t}f\left( s,\Pi _{s}^{n}\right)
ds+\int_{T_{n-1}}^{t}\sigma \left( s,X_{s-}^{n},Y_{s-}^{n}\right) dH_{s}, \\ 
Y_{t}^{n}= & G_{n}\left( X_{T_{n}}^{n}\right) +\int_{t}^{T_{n}}g\left( s,\Pi
_{s}^{n}\right) ds-\int_{t}^{T_{n}}Z_{s}^{n}dH_{s}.%
\end{array}%
\right.  \label{EQ5}
\end{equation}%
Let $L_{G_{n}}$ denotes the Lipschitz constant of the mapping $G_{n}.$ Then,
by Theorem $\ref{STD}$ the required solution of FBSDE $(\ref{EQ5})$ exists
and is unique. Define $G_{n-1}(x)\overset{\vartriangle }{=}Y_{T_{n-1}}^{n},$
then for fixed $x$, $G_{n-1}(x)\in \mathcal{F}_{T_{n-1}}.$ Further, in view
of the Proposition $\ref{Prop2},$ it's straightforward to verify that%
\begin{equation*}
L_{G_{n-1}}\leq \lambda _{1}\overset{\vartriangle }{=}c\left( \left[ \lambda
_{0}+1\right] e^{\left( 2\lambda +\lambda ^{2}\right) \left(
T_{n}-T_{n-1}\right) }-1\right) \leq \bar{\lambda}_{0}.
\end{equation*}%
Next, for $t\in \left[ T_{n-2},T_{n-1}\right] $, we consider the following
FBSDE:%
\begin{equation}
\left\{ 
\begin{array}{ll}
X_{t}^{n-1}= & x+\int_{T_{n-2}}^{t}f\left( s,\Pi _{s}^{n-1}\right)
ds+\int_{T_{n-2}}^{t}\sigma \left( s,X_{s-}^{n-1},Y_{s-}^{n-1}\right) dH_{s},
\\ 
Y_{t}^{n-1}= & G_{n-1}\left( X_{T_{n-1}}^{n-1}\right)
+\int_{t}^{T_{n-1}}g\left( s,\Pi _{s}^{n-1}\right)
ds-\int_{t}^{T_{n-1}}Z_{s}^{n-1}dH_{s}.%
\end{array}%
\right.  \label{EQ6}
\end{equation}%
Once again, since $L_{G_{n-1}}\leq \bar{\lambda}_{0},$ by Theorem $\ref{STD}%
, $ the FBSDE $(\ref{EQ6})$ has a unique solution.

Then as well, we may define $G_{n-2}\left( x\right) ,$ such that%
\begin{equation*}
\begin{array}{r}
L_{G_{n-2}}\leq \lambda _{2}\overset{\vartriangle }{=}c\left( \left[ \lambda
_{1}+1\right] e^{\left( 2\lambda +\lambda ^{2}\right) \left(
T_{n-1}-T_{n-2}\right) }-1\right) \\ 
=c\left( \left[ \lambda _{0}+1\right] e^{\left( 2\lambda +\lambda
^{2}\right) \left( T_{n}-T_{n-2}\right) }-1\right) \leq \bar{\lambda}_{0}.%
\end{array}%
\end{equation*}%
Repeating this procedure backwardly for $i=n,...,1$, we may define $G_{i}$
such that%
\begin{equation*}
L_{G_{i}}\leq \lambda _{n-i}\overset{\vartriangle }{=}c\left( \left[ \lambda
_{0}+1\right] e^{\left( 2\lambda +\lambda ^{2}\right) \left(
T_{n}-T_{i}\right) }-1\right) \leq \bar{\lambda}_{0}.
\end{equation*}%
As a conclusion, one can repeat the above construction and, after a finite
number of steps, we obtain the required unique solution in each subinterval
of the type $\left[ T_{n-i},T_{n-i}\right] $ for $i=0,...,n.$

Now, for $i=1,2,...,n$ and for any $X_{0}\in L^{2}\left( \mathcal{F}%
_{0}\right) $, we construct a solution for the following FBSDE%
\begin{equation*}
\left\{ 
\begin{array}{ll}
X_{t}= & X_{T_{i-1}}+\int_{T_{i-1}}^{t}f\left( s,\Pi _{s}\right)
ds+\int_{T_{i-1}}^{t}\sigma \left( s,X_{s-},Y_{s-}\right) dH_{s}, \\ 
Y_{t}= & G_{i}\left( X_{T_{i}}\right) +\int_{t}^{T_{i}}g\left( s,\Pi
_{s}\right) ds-\int_{t}^{T_{i}}Z_{s}dH_{s}.%
\end{array}%
\right. t\in \left[ t_{i-1},t_{i}\right]
\end{equation*}%
Obviously this provides a solution to the FBSDE $(\ref{EQ1})$. From the
construction and the uniqueness of each step, it is clear that this solution
is unique.

Now, let us prove $\left( ii\right) $. We denote%
\begin{equation*}
V_{t}^{2}=\left\vert f\left( t,0,0,0\right) \right\vert ^{2}+\left\vert
\sigma \left( t,0,0\right) \right\vert ^{2}+\left\vert g\left(
t,0,0,0\right) \right\vert ^{2}.
\end{equation*}%
From Theorem $\ref{STD}$ and by the definition of $G_{i}$, we get%
\begin{equation*}
E\left\{ \left\vert G_{i-1}\left( 0\right) \right\vert ^{2}\right\} \leq
C_{0}E\left\{ \left\vert G_{i}\left( 0\right) \right\vert
^{2}+\int_{T_{i-1}}^{T_{i}}V_{t}^{2}dt\right\} .
\end{equation*}%
By induction one can easily prove that%
\begin{equation*}
\begin{array}{l}
\max\limits_{0\leq i\leq n}E\left\{ \left\vert G_{i}\left( 0\right)
\right\vert ^{2}\right\} \leq C_{0}^{n}E\left\{ \left\vert \varphi \left(
0\right) \right\vert ^{2}+\int_{0}^{T}V_{t}^{2}dt\right\} \\ 
=CE\left\{ \left\vert \varphi \left( 0\right) \right\vert
^{2}+\int_{0}^{T}V_{t}^{2}dt\right\} .%
\end{array}%
\end{equation*}%
Set $n\leq \frac{T}{\delta }+1$\ is a fixed constant depending only on $%
\lambda ,\lambda _{0}$\ and $T$, then so is $C$. Now for $t\in \left[
T_{0},T_{1}\right] $\textbf{,} by using $\left( ii\right) $ of Theorem $\ref%
{STD}$, we get%
\begin{equation*}
\begin{array}{l}
E\left\{ \sup\limits_{0\leq t\leq T}\left\vert X_{t}\right\vert
^{2}+\sup\limits_{0\leq t\leq T}\left\vert Y_{t}\right\vert ^{2}\right\}
\leq CE\left\{ \left\vert X_{0}\right\vert ^{2}+\left\vert G_{1}\left(
0\right) \right\vert ^{2}+\int_{T_{0}}^{T_{1}}V_{t}^{2}dt\right\} \\ 
\leq CE\left\{ \left\vert X_{0}\right\vert ^{2}+\left\vert \varphi \left(
0\right) \right\vert ^{2}+\int_{0}^{T}V_{t}^{2}dt\right\} .%
\end{array}%
\end{equation*}%
Then by induction one can prove%
\begin{equation}
E\left\{ \sup\limits_{0\leq t\leq T}\left\vert X_{t}\right\vert
^{2}+\sup\limits_{0\leq t\leq T}\left\vert Y_{t}\right\vert ^{2}\right\}
\leq CE\left\{ \left\vert X_{0}\right\vert ^{2}+\left\vert \varphi \left(
0\right) \right\vert +\int_{0}^{T}V_{t}^{2}dt\right\} .  \label{inq1}
\end{equation}%
On the other hand, applying Ito's formula to $Y_{t}$ , we obtain%
\begin{equation*}
\begin{array}{l}
E\left\{ \left\vert Y_{0}\right\vert ^{2}+\int_{0}^{T}\left\vert
Z_{t}\right\vert ^{2}dt\right\} =E\left\{ \left\vert Y_{T}\right\vert
^{2}+2\int_{0}^{T}Y_{t}g\left( t,\Pi _{t}\right) dt\right\} \\ 
\leq E\left\{ \left\vert Y_{T}\right\vert ^{2}+C\int_{0}^{T}\left[
\left\vert g\left( t,0,0,0\right) \right\vert ^{2}+\left\vert
X_{t}\right\vert ^{2}+\left\vert Y_{t}\right\vert ^{2}\right] dt+\frac{1}{2}%
\int_{0}^{T}\left\vert Z\right\vert _{t}^{2}dt\right\} .%
\end{array}%
\end{equation*}%
Therefore

\begin{equation}
E\left\{ \int_{0}^{T}\left\vert Z\right\vert _{t}^{2}dt\right\} \leq
CE\left\{ \left\vert X_{0}\right\vert ^{2}+\left\vert \varphi \left(
0\right) \right\vert +\int_{0}^{T}V_{t}^{2}dt\right\} .  \label{inq2}
\end{equation}%
Finally, combining $(\ref{inq1})$ and $(\ref{inq2})$ leads to $\left\Vert
\Pi \right\Vert ^{2}\leq CV_{0}^{2},$ which achieves the proof.\hfill $%
\square $

\subsection{Proof of stability theorem}

\textbf{Proof of Theorem }$\ref{stability of solution}$. For $0\leq
\varepsilon \leq 1,$ let $\Pi ^{\varepsilon }$ be the solution to the
following FBSDE%
\begin{equation*}
\left\{ 
\begin{array}{ll}
X_{t}^{\varepsilon }= & X_{0}+\varepsilon \Delta X_{0}+\int_{0}^{t}\left(
f^{0}\left( s,\Pi _{s}^{\varepsilon }\right) +\varepsilon \Delta f\left(
s,\Pi _{s}^{1}\right) \right) ds \\ 
& +\int_{0}^{t}\left( \sigma ^{0}\left( s,X_{s-}^{\varepsilon
},Y_{s-}^{\varepsilon }\right) +\varepsilon \Delta \sigma \left(
s,X_{s-}^{1},Y_{s-}^{1}\right) \right) dH_{s}; \\ 
Y_{t}^{\varepsilon }= & \left( \varphi ^{0}\left( X_{T}^{\varepsilon
}\right) +\varepsilon \Delta \varphi \left( X_{T}^{1}\right) \right)
+\int_{t}^{T}\left( g^{0}\left( s,\Pi _{s}^{\varepsilon }\right)
+\varepsilon \Delta g\left( s,\Pi _{s}^{1}\right) \right)
ds-\int_{t}^{T}Z_{s}^{\varepsilon }dH_{s}.%
\end{array}%
\right.
\end{equation*}%
and $\triangledown \Pi ^{\varepsilon }$ be the solution of the following
variational linear FBSDE%
\begin{equation*}
\left\{ 
\begin{array}{ll}
\nabla X_{t}^{\varepsilon }= & \Delta X_{0}+\int_{0}^{t}\left(
f_{x}^{0}\left( s,\Pi _{s}^{\varepsilon }\right) \nabla X_{s}^{\varepsilon
}+f_{y}^{0}\left( s,\Pi _{s}^{\varepsilon }\right) \nabla Y_{s}^{\varepsilon
}+f_{z}^{0}\left( s,\Pi _{s}^{\varepsilon }\right) \nabla Z_{s}^{\varepsilon
}+\Delta f\left( s,\Pi _{s}^{1}\right) \right) ds \\ 
& +\int_{0}^{t}\left( \sigma _{x}^{0}\left( s,X_{s-}^{\varepsilon
},Y_{s-}^{\varepsilon }\right) \nabla X_{s}^{\varepsilon }+\sigma
_{y}^{0}\left( s,X_{s-}^{\varepsilon },Y_{s-}^{\varepsilon }\right) \nabla
Y_{s}^{\varepsilon }+\Delta \sigma \left( s,\Pi _{s}^{1}\right) \right)
dH_{s}; \\ 
\nabla Y_{t}^{\varepsilon }= & \varphi _{x}^{0}\left( X_{T}^{\varepsilon
}\right) +\Delta \varphi \left( X_{T}^{1}\right) +\int_{t}^{T}\left(
g_{x}^{0}\left( s,\Pi _{s}^{\varepsilon }\right) \nabla X_{s}^{\varepsilon
}+g_{y}^{0}\left( s,\Pi _{s}^{\varepsilon }\right) \nabla Y_{s}^{\varepsilon
}+g_{z}^{0}\left( s,\Pi _{s}^{\varepsilon }\right) \nabla Z_{s}^{\varepsilon
}+\Delta g\left( s,\Pi _{s}^{1}\right) \right) ds \\ 
& -\int_{t}^{T}\nabla Z_{s}^{\varepsilon }dH_{s};%
\end{array}%
\right.
\end{equation*}%
Then by Theorem \ref{STD}, the above FBSDEs has a unique solution. Moreover,
a simple calculation shows that%
\begin{equation*}
\Delta \Pi _{t}=\int_{0}^{1}\frac{d}{d\varepsilon }\Pi _{t}^{\varepsilon
}d\varepsilon =\int_{0}^{1}\nabla \Pi _{t}^{\varepsilon }d\varepsilon .
\end{equation*}%
since $\left( f^{0},\sigma ^{0},g^{0}\right) $ satisfies $(\ref{H.3})$, by
Lemma \ref{lem3}$,$ we obtain%
\begin{equation*}
\left\Vert \Delta \Pi ^{\varepsilon }\right\Vert ^{2}\leq CE\left\{
\left\vert \Delta X_{0}\right\vert ^{2}+\left\vert \Delta \varphi \left(
X_{T}^{1}\right) \right\vert ^{2}+\int_{0}^{T}\left[ \left\vert \Delta
f\right\vert ^{2}+\left\vert \Delta \sigma \right\vert ^{2}+\left\vert
\Delta g\right\vert ^{2}\right] \left( t,\Pi _{t}^{1}\right) dt\right\} ,
\end{equation*}%
which implies the desired result.\hfill $\square $

\noindent \textbf{Proof of Corollary }$\ref{corol2}$ Using Theorem $\ref%
{stability of solution}$ we have%
\begin{equation*}
\begin{array}{c}
\left\Vert \Pi ^{n}-\Pi ^{n}\right\Vert ^{2}\leq CE\left\{ \left\vert
X_{0}^{n}-X_{0}^{0}\right\vert ^{2}+\left\vert \varphi ^{n}-\varphi
^{0}\right\vert ^{2}\left( X_{T}^{0}\right) \right. \\ 
\left. +\int_{0}^{T}\left[ \left\vert f^{n}-f^{0}\right\vert ^{2}+\left\vert
\sigma ^{n}-\sigma ^{0}\right\vert ^{2}+\left\vert g^{n}-g^{0}\right\vert
^{2}\right] \left( t,\Pi _{t}^{0}\right) dt\right\} .%
\end{array}%
\end{equation*}%
Thus, the desired result follows immediately, by letting $n$ tend to $0$,
and using the dominated convergence theorem.\hfill $\square $

\subsection{Proof of comparison theorem}

\subsubsection{Some auxiliary results}

In order to prove Proposition $\ref{Prop3}$, we need the following two
Lemmas. Let us introduce the following linear FBSDE 
\begin{equation}
\left\{ 
\begin{array}{ll}
X_{t}= & \int_{0}^{t}\left( \bar{a}_{s}^{1}X_{s}+\bar{b}_{s}^{1}\bar{Y}_{s}+%
\bar{c}_{s}^{1}\bar{Z}_{s}\right) ds+\int_{0}^{t}\left( \bar{a}_{s}^{2}X_{s}+%
\bar{b}_{s}^{2}\bar{Y}_{s}\right) dH_{s}, \\ 
\bar{Y}_{t}= & \int_{t}^{T}\left( \bar{a}_{s}^{3}X_{s}+\bar{b}_{s}^{3}\bar{Y}%
_{s}+\bar{c}_{s}^{3}\bar{Z}_{s}\right) ds-\int_{t}^{T}\bar{Z}_{s}dH_{s}.%
\end{array}%
\right.  \label{EQ72}
\end{equation}

Here, $\bar{Y}_{t}\overset{\bigtriangleup }{=}Y_{t}-P_{t}X_{t}$, $\bar{Z}_{t}%
\overset{\bigtriangleup }{=}Z_{t}-P_{t}\left(
a_{t}^{2}X_{t}+b_{t}^{2}Y_{t}\right) -g_{t}X_{t},$ where $P=E\left( P\right)
+\int_{0}^{T}p_{t}dH_{t}$, $P_{t}\overset{\bigtriangleup }{=}E\left(
P\right) +\int_{0}^{t}p_{t}dH_{t};$and%
\begin{equation*}
\left\{ 
\begin{array}{l}
\bar{a}_{t}^{1}\overset{\bigtriangleup }{=}%
a_{t}^{1}+P_{t}b_{t}^{1}+P_{t}a_{t}^{2}c_{t}^{1}+\left\vert P_{t}\right\vert
^{2}b_{t}^{2}c_{t}^{1}+p_{t}c_{t}^{1}; \\ 
\bar{b}_{t}^{1}\overset{\bigtriangleup }{=}%
b_{t}^{1}+P_{t}b_{t}^{2}c_{t}^{1}=b_{t}^{1}; \\ 
\bar{c}_{t}^{1}\overset{\bigtriangleup }{=}c_{t}^{1}; \\ 
\bar{a}_{t}^{2}\overset{\bigtriangleup }{=}a_{t}^{2}+P_{t}b_{t}^{2}; \\ 
\bar{b}_{t}^{2}\overset{\bigtriangleup }{=}b_{t}^{2}; \\ 
\bar{a}_{t}^{3}\overset{\bigtriangleup }{=}%
a_{t}^{3}+p_{t}a_{t}^{2}+P_{t}a_{t}^{1}+\left(
b_{t}^{3}+p_{t}b_{t}^{2}+P_{t}b_{t}^{1}\right) P_{t} \\ 
\text{ \ \ \ \ \ \ }+\left( c_{t}^{3}+P_{t}c_{t}^{1}\right) \left(
p_{t}+P_{t}a_{t}^{2}+\left\vert P_{t}\right\vert ^{2}b_{t}^{2}\right) ; \\ 
\bar{b}_{t}^{3}\overset{\bigtriangleup }{=}%
b_{t}^{3}+p_{t}b_{t}^{2}+P_{t}b_{t}^{1}+P_{t}b_{t}^{2}c_{t}^{3}+\left\vert
P_{t}\right\vert ^{2}b_{t}^{2}c_{t}^{1}; \\ 
\bar{c}_{t}^{3}\overset{\bigtriangleup }{=}c_{t}^{3}+P_{t}c_{t}^{1}.%
\end{array}%
\right.
\end{equation*}

\begin{Lemm}
\label{lem4} Let $\left( X,Y,Z\right) $ be the solution of LFBSDE $\left( 
\text{\ref{EQ7}}\right) ,$ assume $\beta =0$ and $p\leq C$ $.$Then $\left( X,%
\tilde{Y},\tilde{Z}\right) $ is the solution of the linear FBSDE $\left( \ref%
{EQ72}\right) $.
\end{Lemm}

\noindent \textbf{Proof of Lemma} $\ref{lem4}.$ By the definition of $P_{t}$%
, $\bar{Y}_{t}$ and $\bar{Z}_{t}$, we get%
\begin{equation*}
\begin{array}{l}
dX_{t}=\left( a_{t}^{1}X_{t}+b_{t}^{1}\left( \bar{Y}_{t}+P_{t}X_{t}\right)
+c_{t}^{1}\left( \bar{Z}_{t}+P_{t}a_{t}^{2}X_{t}+P_{t}b_{t}^{2}\left( \bar{Y}%
_{t}+P_{t}X_{t}\right) +p_{t}X_{t}\right) \right) dt \\ 
+\left( a_{t}^{2}X_{t}+b_{t}^{2}\left( \bar{Y}_{t}+P_{t}X_{t}\right) \right)
dH_{t} \\ 
=\left( \bar{a}_{t}^{1}X_{t}+\bar{b}_{t}^{1}\bar{Y}_{t}+\bar{c}_{t}^{1}\bar{Z%
}_{t}\right) dt+\left( \bar{a}_{t}^{2}X_{t}+\bar{b}_{t}^{2}\bar{Y}%
_{t}\right) dH_{t},%
\end{array}%
\end{equation*}%
and%
\begin{equation*}
\begin{array}{c}
d\bar{Y}_{t}=-\left( a_{t}^{3}X_{t}+b_{t}^{3}Y_{t}+c_{t}^{3}Z_{t}\right)
dt+Z_{t}dH_{t}-p_{t}\left( a_{t}^{2}X_{t}+b_{t}^{2}Y_{t}\right) dt \\ 
-P_{t}\left( a_{t}^{1}X_{t}+b_{t}^{1}Y_{t}+c_{t}^{1}Z_{t}\right)
dt-P_{t}\left( a_{t}^{2}X_{t}+b_{t}^{2}Y_{t}\right) dH_{t}-p_{t}X_{t}dH_{t}
\\ 
=\bar{Z}_{t}dH_{t}-\left[ \left(
a_{t}^{3}+p_{t}a_{t}^{2}+P_{t}a_{t}^{1}\right) X_{t}+\left(
b_{t}^{3}+p_{t}b_{t}^{2}+G_{t}b_{t}^{1}\right) \left( \bar{Y}%
_{t}+p_{t}X_{t}\right) \right. \\ 
+\left. \left( c_{t}^{3}+P_{t}c_{t}^{1}\right) \left( \bar{Z}_{t}+\left(
p_{t}+P_{t}a_{t}^{2}\right) X_{t}+P_{t}b_{t}^{2}\left( \bar{Y}%
_{t}+P_{t}X_{t}\right) \right) dt\right] \\ 
=-\left( \bar{a}_{t}^{3}X_{t}+\bar{b}_{t}^{3}\bar{Y}_{t}+\bar{c}_{t}^{3}\bar{%
Z}_{t}\right) dt+\bar{Z}_{t}dH_{t},%
\end{array}%
\end{equation*}%
Is easy to prove that $\bar{a}_{t}^{i},\bar{b}_{t}^{i},\bar{c}_{t}^{i}$ are
bounded and still satisfy the assumptions $(\ref{H.3})$. Then this gives the
desired result.

\begin{Lemm}
\label{lem5} Assume $\alpha =0,$\ $c_{t}^{3}=0$, for some integer $m,$\ we
assume $\frac{1}{m}\leq \kappa _{2}\leq m.$ Then there exist small constants 
$\delta $ and $C$ depending on $\lambda $ and $\lambda _{0},$ such that $%
T\leq \delta ,$ and that for some $\varepsilon >0,$%
\begin{equation*}
\left\vert E\left( PX_{t}+\int_{0}^{T}\left(
a_{t}^{3}X_{t}+b_{t}^{3}Y_{t}\right) dt\right) \right\vert \leq Cm\sqrt{%
\varepsilon }T.
\end{equation*}
\end{Lemm}

\noindent \textbf{Proof of Lemma} $\ref{lem5}.$ By standard arguments and
using Young's inequality, for every $\varepsilon >0$, there exist constant $%
C $\ depending only on $\lambda ,\lambda _{0,}$ that%
\begin{equation*}
\begin{array}{l}
\sup\limits_{0\leq t\leq T}E\left( \left\vert X_{t}\right\vert
^{2}+\left\vert Y_{t}\right\vert ^{2}\right) +E\left( \int_{0}^{T}\left\Vert
Z_{t}\right\Vert _{l^{2}\left( \mathbb{R}\right) }^{2}dt\right) \leq
C\varepsilon ^{-1}E\left( \int_{0}^{T}\left( \left\vert X_{t}\right\vert
^{2}+\left\vert Y_{t}\right\vert ^{2}\right) dt\right) +\frac{\varepsilon }{2%
}E\left( \int_{0}^{T}\left\vert \beta _{t}\right\vert ^{2}dt\right) \\ 
\leq C\varepsilon ^{-1}T\sup_{0\leq t\leq T}E\left( \left\vert
X_{t}\right\vert ^{2}+\left\vert Y_{t}\right\vert ^{2}\right) +\frac{%
\varepsilon }{2}m^{2}T.%
\end{array}%
\end{equation*}%
If we choose the constant $\delta =\frac{\epsilon }{2C}$\ and will specify $%
\varepsilon $\ later. Then for $T\leq \delta $, we get%
\begin{equation*}
\sup_{0\leq t\leq T}E\left( \left\vert X_{t}\right\vert ^{2}+\left\vert
Y_{t}\right\vert ^{2}\right) +E\left( \int_{0}^{T}\left\Vert
Z_{t}\right\Vert _{l^{2}\left( \mathbb{R}\right) }^{2}dt\right) \leq
m^{2}\varepsilon T.
\end{equation*}%
And%
\begin{equation*}
\begin{array}{l}
E\left( \left\vert X_{t}\right\vert ^{2}\right) \leq CE\left( \left\vert
\int_{0}^{T}\left( a_{t}^{1}X_{t}+b_{t}^{1}Y_{t}+c_{t}^{1}Z_{t}\right)
dt\right\vert ^{2}+\left\vert \int_{0}^{T}\left(
a_{t}^{2}X_{t}+b_{t}^{2}Y_{t}\right) dH_{t}\right\vert ^{2}\right) \\ 
\leq CE\left( T\int_{0}^{T}\left( \left\vert X_{t}\right\vert
^{2}+\left\vert Y_{t}\right\vert ^{2}+\left\Vert Z_{t}\right\Vert
_{l^{2}\left( \mathbb{R}\right) }^{2}\right) dt+\int_{0}^{T}\left(
\left\vert X_{t}\right\vert ^{2}+\left\vert Y_{t}\right\vert ^{2}\right)
dt\right) \\ 
\leq Cm^{2}\varepsilon T^{2}.%
\end{array}%
\end{equation*}%
Thus%
\begin{equation*}
\begin{array}{l}
\left\vert E\left( PX_{t}\right) +\int_{0}^{T}\left(
a_{t}^{3}X_{t}+b_{t}^{3}Y_{t}\right) dt\right\vert \\ 
\leq CE^{\frac{1}{2}}\left( \left\vert X_{T}\right\vert ^{2}\right)
+CT\sup\limits_{0\leq t\leq T}E^{\frac{1}{2}}\left( \left\vert
X_{t}\right\vert ^{2}+\left\vert Y_{t}\right\vert ^{2}\right) \leq Cm\sqrt{%
\varepsilon }T.%
\end{array}%
\end{equation*}%
This ends the proof.\hfill $\square $

\subsubsection{\textbf{Proof of Proposition.}\protect\ref{Prop3}\textbf{. }}

The proof of the proposition $\ref{Prop3}$ will be splitted into several
steps.

$Step$ $1.$ Assume that $P=0$ and $\beta =0$. If $Y_{0}<0,$ let us define
the following stopping time%
\begin{equation*}
\tau \overset{\bigtriangleup }{=}\inf \left\{ t:Y_{t}=0\right\} \wedge T.
\end{equation*}%
Since $Y_{T}=\alpha \geq 0,$ we get $Y\tau =0.$ Define%
\begin{equation*}
\hat{a}_{t}^{i}\overset{\bigtriangleup }{=}a_{t}^{i}1_{\left\{ \tau
>t\right\} };\hat{b}_{t}^{i}\overset{\bigtriangleup }{=}b_{t}^{i}1_{\left\{
\tau >t\right\} };\hat{c}_{t}^{i}\overset{\bigtriangleup }{=}%
c_{t}^{i}1_{\left\{ \tau >t\right\} }
\end{equation*}%
\begin{equation*}
\hat{X}_{t}\overset{\bigtriangleup }{=}X_{\tau \wedge t};\text{ \ \ }\hat{Y}%
_{t}\overset{\bigtriangleup }{=}Y_{\tau \wedge t};\text{ \ \ }\hat{Z}_{t}%
\overset{\bigtriangleup }{=}Z_{\tau \wedge t}\text{\ }
\end{equation*}%
In view of Lemma $\ref{lem3}$, the following LFBSDE:%
\begin{equation*}
\left\{ 
\begin{array}{ll}
\hat{X}_{t}= & \int_{0}^{t}\left( \hat{a}_{s}^{1}\hat{X}_{s}+\hat{b}_{s}^{1}%
\hat{Y}_{s}+\hat{c}_{s}^{1}\hat{Z}_{s}\right) ds+\int_{0}^{t}\left( \hat{a}%
_{s}^{2}\hat{X}_{s}+\hat{b}_{s}^{2}\hat{Y}_{s}\right) dH_{s}, \\ 
\hat{Y}_{t}= & \int_{t}^{T}\left( \hat{a}_{s}^{3}\hat{X}_{s}+\hat{b}_{s}^{3}%
\hat{Y}_{s}+\hat{c}_{s}^{3}\hat{Z}_{s}\right) ds-\int_{t}^{T}\hat{Z}%
_{s}dH_{s},%
\end{array}%
\right.
\end{equation*}%
has a unique solution, with $\hat{Y}_{T}=0.$That is to say $Y_{0}=\hat{Y}%
_{0}=0,$ obviously this leads to a contradiction. In other words, we have
proved that $Y_{0}\geq 0.$\medskip

$Step$ $2.$ Assume that all the conditions in Lemma $\ref{lem4}$ are
fulfilled, then $\bar{Y}_{T}=\alpha \geq 0$. Applying $Step$ $1$ we get $%
Y_{0}=\hat{Y}_{0}\geq 0.$\medskip

$Step$ $3.$ Assume $\beta =0.$ One can find $P_{n}$ satisfying the condition
in Lemma $\ref{lem4}$ such that $P_{n}\rightarrow P$ a.s. and $\left\vert
P_{n}\right\vert \leq \lambda .$ Let $\left( X^{n},Y^{n},Z^{n}\right) $
denotes the solution corresponding to $G_{n}.$ Apply the result of $Step2$
to conclude that $Y_{0}^{n}\geq 0$. Then from Corollary $\ref{corol2},$ we
get $Y_{0}=\lim\limits_{n\rightarrow \infty }Y_{0}^{n}\geq 0.$\medskip

$Step$ $4.$ Assume all the conditions in Lemma $\ref{lem5}$ are in force.\
Then%
\begin{equation*}
\begin{array}{l}
Y_{0}=E\left( PX_{T}+\int_{0}^{T}\left( a_{t}^{3}X_{t}+b_{t}^{3}Y_{t}+\beta
_{t}\right) dt\right) \\ 
\geq m^{-1}T-\left\vert E\left( PX_{T}+\int_{0}^{T}\left(
a_{t}^{3}X_{t}+b_{t}^{3}Y_{t}\right) dt\right) \right\vert \\ 
\geq m^{-1}T-Cm\sqrt{\varepsilon }T.%
\end{array}%
\end{equation*}%
Now choose $\varepsilon =C^{-2}m^{-4}$, we get $Y_{0}\geq 0$.\medskip

$Step$ $5.$ Assume $\frac{1}{m}\leq \beta \leq m$ and $T\leq \delta ,$\
where $\delta $ is the same as in Lemma $\ref{lem5}$. Denote%
\begin{equation*}
\left\{ 
\begin{array}{cc}
X_{t}^{\prime }= & \int_{0}^{t}\left( a_{s}^{1}X_{s}^{\prime
}+b_{s}^{1}Y_{s}^{\prime }+c_{s}^{1}Z_{s}^{\prime }\right)
ds+\int_{0}^{t}\left( a_{s}^{2}X_{s}^{\prime }+b_{s}^{2}Y_{s}^{\prime
}\right) dH_{s}, \\ 
Y_{t}^{\prime }= & PX_{T}^{\prime }+\alpha +\int_{t}^{T}\left(
a_{s}^{3}X_{s}^{\prime }+b_{s}^{3}Y_{s}^{\prime }+c_{s}^{3}Z_{s}^{\prime
}\right) ds-\int_{t}^{T}Z_{s}^{\prime }dH_{s},%
\end{array}%
\right.
\end{equation*}%
and%
\begin{equation*}
\left\{ 
\begin{array}{cc}
X_{t}^{\prime \prime }= & \int_{0}^{t}\left( a_{s}^{1}X_{s}^{\prime \prime
}+b_{s}^{1}Y_{s}^{\prime \prime }+c_{s}^{1}Z_{s}^{\prime \prime }\right)
ds+\int_{0}^{t}\left( a_{s}^{2}X_{s}^{\prime \prime }+b_{s}^{2}Y_{s}^{\prime
\prime }\right) dH_{s}, \\ 
Y_{t}^{\prime \prime }= & LX_{T}^{\prime \prime }+\int_{t}^{T}\left(
a_{s}^{3}X_{s}^{\prime \prime }+b_{s}^{3}Y_{s}^{\prime \prime
}+c_{s}^{3}Z_{s}^{\prime \prime }+\beta _{s}\right)
ds-\int_{t}^{T}Z_{s}^{\prime \prime }dH_{s},%
\end{array}%
\right.
\end{equation*}%
By Step $3$, $Y_{0}^{\prime }\geq 0$ , and by Step $4$, $Y_{0}^{\prime
\prime }\geq 0$. Then, $Y_{0}=Y_{0}^{\prime }+Y_{0}^{\prime \prime }\geq 0.$%
\medskip

$Step$ $6.$ Assume $\frac{1}{m}\leq \beta \leq m.$ Let $\delta $ be as in
Lemma $\ref{lem5}$ but corresponding to $\left( \lambda ,\bar{\lambda}%
_{0},m\right) $ instead of $\left( \lambda ,\lambda _{0},m\right) $, and
assume $\left( n-1\right) \delta <T<n\delta $. Denote $T_{i}\overset{%
\bigtriangleup }{=}\frac{iT}{n}$ , $L_{n}\overset{\bigtriangleup }{=}L$ and $%
\alpha _{n}\overset{\bigtriangleup }{=}\alpha .$ For $t\in \left[
T_{n-1},T_{n}\right] $, let%
\begin{equation*}
\left\{ 
\begin{array}{ll}
X_{t}^{n,1}= & 1+\int_{T_{n-1}}^{t}\left(
a_{s}^{1}X_{s}^{n,1}+b_{s}^{1}Y_{s}^{n,1}+c_{s}^{1}Z_{s}^{n,1}\right)
ds+\int_{T_{n-1}}^{t}\left( a_{s}^{2}X_{s}^{n,1}+b_{s}^{2}Y_{s}^{n,1}\right)
dH_{s}, \\ 
Y_{t}^{n,1}= & P_{n}X_{T}^{n,1}+\int_{t}^{T_{n}}\left(
a_{s}^{3}X_{s}^{n,1}+b_{s}^{3}Y_{s}^{n,1}+c_{s}^{3}Z_{s}^{n,1}\right)
ds-\int_{t}^{T_{n}}Z_{s}^{n,1}dH_{s},%
\end{array}%
\right.
\end{equation*}%
and%
\begin{equation*}
\left\{ 
\begin{array}{cc}
X_{t}^{n,2}= & 1+\int_{T_{m-1}}^{t}\left(
a_{s}^{1}X_{s}^{n,2}+b_{s}^{1}Y_{s}^{n,2}+c_{s}^{1}Z_{s}^{n,2}\right)
ds+\int_{T_{n-1}}^{t}\left( a_{s}^{2}X_{s}^{n,2}+b_{s}^{2}Y_{s}^{n,2}\right)
dH_{s}, \\ 
Y_{t}^{n,2}= & P_{n}X_{T}^{n}+\alpha _{n}+\int_{t}^{T_{n}}\left(
a_{s}^{3}X_{s}^{n,2}+b_{s}^{3}Y_{s}^{n,2}+c_{s}^{3}Z_{s}^{n,2}+\beta
_{s}\right) ds-\int_{t}^{T_{n}}Z_{s}^{n,2}dH_{s},%
\end{array}%
\right.
\end{equation*}%
Denote 
\begin{equation*}
P_{n-1}\overset{\bigtriangleup }{=}Y_{T_{n-1}}^{n,1},\alpha _{n-1}\overset{%
\bigtriangleup }{=}Y_{T_{n-1}}^{n,2}.
\end{equation*}%
By the proof of Theorem $\ref{LTD},$ we know that $\left\vert
P_{n-1}\right\vert \leq \lambda _{1}\leq \bar{\lambda}_{0}.$ Apply the
result of $Step$ $5$, we get $\alpha _{n-1}\geq 0.$ We note that, for\ $t\in %
\left[ 0,T_{n-1}\right] $, $\left( X,Y,Z\right) $ satisfies%
\begin{equation*}
\left\{ 
\begin{array}{ll}
X_{t}= & \int_{0}^{t}\left(
a_{s}^{1}X_{s}+b_{s}^{1}Y_{s}+c_{s}^{1}Z_{s}\right) ds+\int_{0}^{t}\left(
a_{s}^{2}X_{s}+b_{s}^{2}Y_{s}\right) dH_{s}, \\ 
Y_{t}= & P_{n-1}X_{T_{n-1}}+\alpha _{n-1}+\int_{t}^{T_{n-1}}\left(
a_{s}^{3}X_{s}+b_{s}^{3}Y_{s}+c_{s}^{3}Z_{s}+\beta _{s}\right)
ds-\int_{t}^{T_{n-1}}Z_{s}dH_{s}.%
\end{array}%
\right.
\end{equation*}%
Repeating the same arguments, we may define $L_{1}$ and $\alpha _{1}\geq 0,$
and it holds that%
\begin{equation*}
\left\{ 
\begin{array}{ll}
X_{t}= & \int_{0}^{t}\left(
a_{s}^{1}X_{s}+b_{s}^{1}Y_{s}+c_{s}^{1}Z_{s}\right) ds+\int_{0}^{t}\left(
a_{s}^{2}X_{s}+b_{s}^{2}Y_{s}\right) dH_{s}, \\ 
Y_{t}= & P_{1}X_{T_{1}}+\alpha _{1}+\int_{t}^{T_{1}}\left(
a_{s}^{3}X_{s}+b_{s}^{3}Y_{s}+c_{s}^{3}Z_{s}+\beta _{s}\right)
ds-\int_{t}^{T_{1}}Z_{s}dH_{s}.%
\end{array}%
\right.
\end{equation*}%
By step $5$, we have $Y_{0}\geq 0.$\medskip

$Step$ $7$. In the general case, we put $\beta ^{m}\overset{\bigtriangleup }{%
=}\left( \beta \wedge m\right) \vee \frac{1}{m}$ and let $\left(
X^{m},Y^{m},Z^{m}\right) $ denote the solution corresponding to $\beta ^{m}.$
We know by Step $6,$ that$\ Y_{0}^{m}\geq 0.$ Then by Corollary $\ref{corol2}%
,$ $Y_{0}=\lim\limits_{m\rightarrow \infty }Y_{0}^{m}\geq 0.$ This gives the
result.\hfill $\square $\bigskip

We are now in position to give the proof of comparison theorem.

\subsubsection{\textbf{Proof of Theorem }\protect\ref{comparison of
solutions}\textbf{. }}

For $0\leq \varepsilon \leq 1,$ let $\Pi ^{\varepsilon }$ and $\nabla \Pi
^{\varepsilon \text{ }}$ be as in the prove of Theorem $\ref{stability of
solution}.$ Then, we get $\Delta X_{0}=0,\Delta f=0,\Delta \sigma =0,\Delta
g\geq 0,\ \Delta \varphi \geq 0.$\ From Proposition $\ref{Prop3}$, we have $%
\nabla Y_{0}^{\varepsilon }\geq 0.$ This proves the theorem.\hfill $\square $%


\bigskip

\begin{thebibliography}{99}
\bibitem{A} F. Antonelli, \textit{Backward-forward stochastic differential
equations.} Ann. Appl. Probab. 3 $\left( 1993\right) $, 777-793.

\bibitem{B} F. Baghery, N. Khelfallah, B. Mezerdi and I. Turpin, Fully
coupled \textit{forward backward stochastic differential equations driven by
Lévy processes and application to differential games}, Rand. Operat. Stoch.
Equa. 22 $\left( 2014\right) $, 151--161.

\bibitem{BE} K. Bahlali, M. Eddahbi and E. Essaky, \textit{BSDE associated
with Lévy processes and application to PDIE}, Journal of Applied Mathematics
and Stochastic Analysis. 16, part 1 $(2003)$, 1--17.

\bibitem{D} F. Delarue, \textit{On the existence and uniqueness of solutions
to FBSDEs in a non-degenerate case}, Stoch. Proc. Appl. 99 $\left(
2002\right) $, 209--286.

\bibitem{H} S. Hamadène, \textit{Backward-Forward SDE's and stochastic
differential game}, Stoch. Proc. Appl. 77 $\left( 1998\right) $, 1-15.

\bibitem{HP} Y. Hu and S. Peng, \textit{Solution of forward-backward
stochastic differential equations}, Probab. Th. Rel. Fields. 103 $\left(
1995\right) $, 273--283.

\bibitem{MP} J. Ma, P. Protter and J.\ Yong, \textit{Solving
forward-backward stochastic differential equations explicitly - a four step
scheme,} Probab. Th. Rel. Fields. 98 $\left( 1994\right) $, 339--359.

\bibitem{NS} D. Nualart and W. Schoutens, \textit{Chaotic and predictable
representation for Lévy processes with applications in finance}, Stoch..
Process. Appl. 90 $\left( 2000\right) $, 109--122.

\bibitem{NSS} D. Nualart and W. Schoutens, \textit{BSDEs and Feynman-Kac
formula for Lévy processes with applications in finance}, Bernoulli. 7 $%
\left( 2001\right) $, 761--776.

\bibitem{PT} E. Pardoux and S. Tang, \textit{Forward-backward stochastic
differential equations and quasilinear parabolic PDEs, }Probab. Theory
Related Fields. 114 $\left( 1999\right) $, 123--150.

\bibitem{PW} S. Peng and Z.\ Wu, \textit{Fully coupled forward-backward
stochastic differential equations and applications to optimal control}, SIAM
J. Control Optim. 37 $\left( 1999\right) $, 825--843.

\bibitem{PS} R.S. Pereire and E.\ Shamarova, \textit{Forward-backward SDEs
driven by Lévy processes and application to option pricing, }Global and
Stochastic Analysis. 2 $\left( 2012\right) $, 2248--9444.\textit{\ }

\bibitem{P} P. Protter, \textit{Stochastic Integration and Differential
Equations,} Spring-Verlag, Berlin. $\left( 1990\right) $.

\bibitem{WZ} Z. Wu, \textit{Fully coupled FBSDE with Brownian motion and
Poisson process in stopping time duration.} J. Austr. Math. Soc. (2003),
74.02: 249-266.

\bibitem{Y} J.\ Yong, \textit{Finding adapted solutions of forward-backward
stochastic differential equations - method of continuation}, Probab. Theory
Related Fields. 107 $\left( 1997\right) $, 537--572.

\bibitem{Z} J. Zhang, \textit{The Wellposedness of FBSDEs. }Discrete and
continuous, Dynamical System.$\left( 2006\right) $, 927-940.
\end{thebibliography}
\end{document}